\documentclass[oneside,11pt]{article}

\usepackage{epsfig}
\usepackage{natbib}
\usepackage{amsfonts}
\usepackage{amsmath}
\usepackage{amssymb}
\usepackage{bbm}
\usepackage{multirow}
\usepackage{longtable}
\usepackage{tabularx}
\usepackage{graphicx}
\newcommand{\inner}[1]{\left\langle#1\right\rangle}
\def\H{\mathcal{H}}
\def\HS{\mathbb{H}}
\def\E{\mathcal{E}}

\def\A{\mathcal{A}}

\def\R{\mathbb{R}}

\def\N{\mathbb{N}}

\def\B{\mathcal{B}}

\newcommand{\norm}[1]{\left\|#1\right\|}

\def\bydef{:=}
\def\Pr{\mathrm{P}}

\def\det{\mathop{\rm det}\nolimits}

\def\Exp{\mathbb{E}\mathop{\!}\nolimits}
\def\vol{\mathop{\rm vol}\nolimits}
\def\div{\mathop{\rm div}\nolimits}
\def\grad{\mathop{\rm grad}\nolimits}
\def\inj{\mathrm{inj}}
\newcommand{\Id}{\mathbbm{1}}

\def\Delrw{\Delta^{\text{(rw)}}}
\def\Delu{\Delta^{\text{(u)}}}
\def\Deln{\Delta^{\text{(n)}}}

\def\eps{\varepsilon}

\def\Var{\mathop{\rm Var}\nolimits}%

\newtheorem{theorem}{Theorem}
\newtheorem{lemma}[theorem]{Lemma}

\newenvironment{remark}{\par\noindent{\bf Remark:\ }}{}
\newtheorem{definition}[theorem]{Definition}
\newtheorem{corollary}[theorem]{Corollary}
\newtheorem{proposition}[theorem]{Proposition}
\newtheorem{assumptions}[theorem]{Assumption}
\newenvironment{proof}{\par\noindent{\bf Proof:\ }}{\hfill$\Box$\\[2mm]}

\def\dsE{\mathbb{E}}

\def\begar{$$\begin{array}{lll}}
\def\endar{\end{array}$$}
\def\begarlab{\begin{equation} \begin{array}{lll} \label}
\def\endarlab{\end{array} \end{equation}}
\def\lbegar{$$\left\{ \begin{array}{lll}}
\def\rendar{\end{array} \right.$$}

\newcommand\lam{\lambda}

\def\td{\tilde{d}}
\def\tk{\tilde{k}}
\def\tA{\tilde{A}}

\newlength{\myVSpace}
\newcommand\xstrut{\raisebox{-.25\myVSpace}
  {\rule{0pt}{\myVSpace}}%
}







\makeatletter
\newenvironment{small-lar}[1]{\skip@=\baselineskip#1 \baselineskip=\skip@\lbegar }{\rendar \ignorespacesafterend}
\makeatother

\makeatletter
\newenvironment{small-ar}[1]{\skip@=\baselineskip#1 \baselineskip=\skip@\begar }{\endar \ignorespacesafterend}
\makeatother

\begin{document}

\title{Graph Laplacians and their convergence on random neighborhood graphs}

\author{Matthias Hein,\\ Max Planck Institute for Biological Cybernetics, T{\"u}bingen, Germany \and
   Jean-Yves Audibert,\\  CERTIS, ENPC, Paris, France  \and
   Ulrike von Luxburg,\\ Max Planck Institute for Biological Cybernetics, T{\"u}bingen, Germany}


\maketitle

\begin{abstract}
Given a sample from a probability measure with support on a
submanifold in Euclidean space one can construct a neighborhood
graph which can be seen as an approximation of the submanifold.
The graph Laplacian of such a graph is used in several machine
learning methods like semi-supervised learning, dimensionality
reduction and clustering. In this paper we determine the pointwise
limit of three different graph Laplacians used in the literature
as the sample size increases and the neighborhood size approaches
zero. We show that for a uniform measure on the submanifold all
graph Laplacians have the same limit up to constants. However in
the case of a non-uniform measure on the submanifold only the so
called random walk graph Laplacian converges to the weighted
Laplace-Beltrami operator.
\end{abstract}

\section{Introduction}
In recent years, methods based on graph Laplacians have become
increasingly popular in machine learning. They have been used in
semi-supervised learning \citep{BelNiy03b,Zhou2003,ZhuGha2002},
spectral clustering \citep{SpiTen1996,Luxburg06_TR} and
dimensionality reduction \citep{BelNiy03,CoiLaf2005}. Their
popularity is mainly due to the following properties of the
Laplacian which will be discussed in more detail in a later
section:
\begin{itemize}
\item the Laplacian is the generator of the diffusion process
(label propagation in semi-supervised learning), \item the
eigenvectors of the Laplacian have special geometric properties
(motivation for spectral clustering), \item the Laplacian induces
an adaptive regularization functional, which adapts to the density
and the geometric structure
      of the data (semi-supervised learning, classification).
\end{itemize}

If the data lies in $\R^d$ the neighborhood graph built from the
random sample can be seen as an approximation of the continuous
structure. in particular,  if the data has support on a
low-dimensional submanifold the neighborhood graph is a discrete
approximation of the submanifold. In machine learning we are
interested in the intrinsic properties and objects of this
submanifold. The approximation of the Laplace-Beltrami operator
via the graph Laplacian is a very important one since it has
numerous applications as we will discuss later.

Approximations of the Laplace-Beltrami operator or related objects
have been studied for certain special deterministic graphs. The
easiest case is a grid in $\R^d$. In numerics it is standard to
approximate the Laplacian with finite-differences schemes on the
grid. These can be seen as a special instances of a graph
Laplacian. There convergence for decreasing grid-size follows
easily by an argument using Taylor expansions. Another more
involved example is the work of \cite{Var1984}, where for a graph
generated by an $\epsilon$-packing of a manifold, the equivalence
of certain properties of random walks on the graph and Brownian
motion on the manifold have been established. The connection
between random walks and the graph Laplacian becomes obvious by
noting that the graph Laplacian as well as the Laplace-Beltrami
operator are the generators of the diffusion process on the graph
and the manifold, respectively. In \cite{Xu2004} the convergence
of a discrete approximation of the Laplace Beltrami operator for a
triangulation of a 2D-surface in $\R^3$ was shown. However,  it is
unclear whether the approximation described there can be written
as a graph Laplacian and whether this result can be generalized to
higher dimensions.

In the case where the graph is generated randomly, only first
results have been proved so far. The first work on the large
sample limit of graph Laplacians has been done by
\cite{BouChaHei2003}. There the authors studied the convergence of
the regularization functional induced by the graph Laplacian using
the law of large numbers for $U$-statistics. In a second step
taking the limit of the neighborhoodsize $h \rightarrow 0$, they
got $\frac{1}{p^2}\nabla (p^2\nabla)$ as the effective limit
operator in $\R^d$. Their result has recently been generalized to
the submanifold case and uniform convergence over the space of
H\"older-functions by \cite{Hei2005,Hei2006}. In
\cite{LuxBelBou2004}, the neighborhoodsize $h$ was kept fixed
while the large sample limit of the graph Laplacian was
considered. In this setting, the authors showed strong convergence
results of graph Laplacians to certain integral operators, which
imply the convergence of the eigenvalues and eigenfunctions.
Thereby showing the consistency of spectral clustering for a fixed
neighborhood size.

In contrast to the previous work  in this paper we will consider
the large sample limit and the limit as the neighborhood size
approaches zero simultaneously for a certain class of
neighbhorhood graphs. The main emphasis lies on the case where the
data generating measure has support on a submanifold of $\R^d$.
The bias term, that is the difference between the continuous
counterpart of the graph Laplacian and the Laplacian itself has
been studied first for compact submanifolds without boundary by
\cite{SmoWeiWit2000} and \cite{Belkin2003} for the Gaussian kernel
and a uniform data generating measure and was then generalized by
\cite{Lafon2004} to general isotropic weights and general
probability measures. Additionally Lafon showed that the use of
data-dependent weights for the graph allows to control the
influence of the density. They all show that the bias term
converges pointwise if the neighborhood size goes to zero. The
convergence of the graph Laplacian towards these continuous
averaging operators was left open. This part was first studied by
\cite{HeiAudLux2005} and \cite{BelNiy2005}. In \cite{BelNiy2005}
the convergence was shown for the so called unnormalized graph
Laplacian in the case of a uniform probability measure on a
compact manifold without boundary and using the Gaussian kernel
for the weights, whereas in \cite{HeiAudLux2005} the pointwise
convergence was shown for the random walk graph Laplacian in the
case of general probability measures on non-compact manifolds with
boundary using general isotropic data-dependent weights. More
recently \cite{GinKol2006} have extended the pointwise convergence
for the unnormalized graph Laplacian shown by \cite{BelNiy2005} to
uniform convergence on compact submanifolds without boundary
giving explicit rates. In \cite{Sin2006}, see also
\cite{GinKol2006}, the rate of convergence given by
\cite{HeiAudLux2005} has been improved in the setting of the
uniform measure. In this paper we will study the three most often
used graph Laplacians in the machine learning literature and show
their pointwise convergence in the general setting of
\cite{Lafon2004} and \cite{HeiAudLux2005}, that is we will in
particular consider the case where by using data-dependent weights
for the graph we can control the influence of the density on the
limit operator.

In Section \ref{sec:graphs} we introduce the basic framework
necessary to define graph Laplacians for general directed weighted
graphs and then simplify the general case to undirected graphs. in
particular,  we define the three graph Laplacians used in machine
learning so far, which we call the normalized, the unnormalized
and the random walk Laplacian. In Section \ref{sec:limit} we
introduce the neighborhood graphs studied in this paper, followed
by an introduction to the so called weighted Laplace-Beltrami
operator, which will turn out to be the limit operator in general.
We also study properties of this limit operator and provide
insights why and how this operator can be used for semi-supervised
learning, clustering and regression. Then finally we present the
main convergence result for all three graph Laplacians and give
the conditions on the neighborhood size as a function of the
sample size necessary for convergence. In Section
\ref{sec:illustration} we illustrate the main result by studying
the difference between the three graph Laplacians and the effects
of different data-dependent weights on the limit operator. In
Section \ref{sec:continuum-limit} we prove the main result. We
introduce a framework for studying non-compact manifolds with
boundary and provide the necessary assumptions on the submanifold
$M$, the data generating measure $P$ and the kernel $k$ used for
defining the weights of the edges. We would like to note that the
theorems given in Section \ref{sec:continuum-limit} contain
slightly stronger results than the ones presented in Section
\ref{sec:limit}. The reader who is not familiar with differential
geometry will find a brief introduction to the basic material used
in this paper in Appendix \ref{sec:diffgeometry}.

\section{Abstract Definition of the Graph Structure}\label{sec:graphs}
In this section we define the structure on a graph which is
required in order to define the graph Laplacian. To this end one
has to introduce Hilbert spaces $H_V$ and $H_E$ of functions on
the vertices $V$ and edges $E$, define a difference operator $d$,
and then set the graph Laplacian as $\Delta=d^*d$. We first do
this in full generality for directed graphs and then specialize it
to undirected graphs. This approach is well-known for undirected
graphs in discrete potential theory and spectral graph theory, see
for example \cite{Dod1984,Chung1997,Woe2000,DonMey2002}, and was
generalized to directed graphs by \cite{ZhoSchHof2005}
for a special choice of $H_V,H_E$ and $d$. To our knowledge the very general setting introduced here has not been discussed elsewhere.\\
In many articles graph Laplacians are used without explicitly
mentioning $d$, $H_V$ and $H_E$. This can be misleading since, as
we will show, there always exists a whole family of choices for
$d$, $H_V$ and $H_E$ which all yield the same graph Laplacian.

\subsection{Hilbert Spaces of Functions on the Vertices $V$ and the Edges $E$}\label{subsec:Hilbert spaces graphs}
Let $(V,W)$ be a graph where $V$ denotes the set of vertices with
$|V|=n$ and $W$ a positive $n \times n$ similarity matrix, that is
$w_{ij}\geq 0, \; i,j=1,\ldots,n$. $W$ need not to be symmetric,
that means we consider the case of a directed graph. The special
case of an undirected graph will be discussed in a following
section. Let $E \subset V \times V$ be the set of edges
$e_{ij}=(i,j)$ with $w_{ij}>0$. $e_{ij}$ is said to be a directed
edge from the vertex $i$ to the vertex $j$ with weight $w_{ij}$.
Moreover,  we define the outgoing and ingoing sum of weights of a
vertex $i$ as
\begin{align*}
 d^{out}_i=\frac{1}{n}\sum_{j=1}^n w_{ij}, \quad d^{in}_i=\frac{1}{n}\sum_{j=1}^n w_{ji}.
\end{align*}
We assume that $d^{out}_i+d^{in}_i>0, \; i=1,\ldots,n$, meaning
that each vertex has at least one in- or outgoing edge. Let
$\R_+=\{x \in \R \, | \, x\geq 0\}$ and $\R_+^*=\{x \in \R \, | \,
x>0\}$. The inner product on the function space $\R^V$ is defined
as
\[ \inner{f,g}_V=\frac{1}{n}\sum_{i=1}^n\, f_i \,g_i\, \chi_i,\]
where $\chi_i=(\chi_{out}(d^{out}_i) + \chi_{in}(d^{in}_i))$ with
$\chi_{out}:\R_+ \rightarrow \R_+$ and $\chi_{in}:\R_+ \rightarrow
\R_+$,  $\chi_{out}(0)=\chi_{in}(0)=0$ and further $\chi_{out}$ and $\chi_{in}$ strictly positive on $\R_+^*$.\\
We also define an inner product on the space of functions $\R^E$
on the edges:
\[ \inner{F,G}_E = \frac{1}{2n^2}\sum_{i,j=1}^n \, F_{ij}\, G_{ij} \,\phi(w_{ij}),\]
where $\phi:\R_+ \rightarrow \R_+$, $\phi(0)=0$ and $\phi$
strictly positive on $\R_+^*$. Note that with these assumptions on
$\phi$ the sum is taken only over the set of edges $E$. One can
check that both inner products are well-defined. We denote by
$\H(V,\chi)=(\R_V,\inner{\cdot,\cdot}_V)$ and
$\H(E,\phi)=(\R^E,\inner{\cdot,\cdot}_E)$ the corresponding
Hilbert spaces. As a last remark let us clarify the roles of
$\R^V$ and $\R^E$. The first one is the space of functions on the
vertices and therefore can be regarded as a normal function space.
However,  elements of $\R^E$ can be interpreted as a flow on the
edges so that the function value on an edge $e_{ij}$ corresponds
to the "mass" flowing from one vertex $i$ to the vertex $j$ (per
unit time).
\subsection{The Difference Operator $d$ and its Adjoint $d^*$}\label{subsec:difference operator}
\begin{definition} \label{def:difference}
The \textbf{difference operator} $d: \H(V,\chi) \rightarrow
\H(E,\phi)$ is defined as follows:
\begin{equation}
\forall \; e_{ij} \in E, \quad \quad
(df)(e_{ij})=\gamma(w_{ij})\,(f(j)-f(i)),
\end{equation}
where $\gamma:\R^*_+ \rightarrow \R^*_+$.
\end{definition}
\begin{remark}
Note that $d$ is zero on the constant functions as one would
expect it from a derivative. In \cite{Zhou2003} another difference
operator $d$ is used:
\begin{equation}\label{eq:denny-difference}
 (df)(e_{ij})=\gamma(w_{ij})\Big( \tfrac{f(j)}{\sqrt{d(j)}}-\tfrac{f(i)}{\sqrt{d(i)}} \Big).
\end{equation}
Note that in \cite{Zhou2003} they have $\gamma(w_{ij})\equiv 1$.
This difference operator is in general not zero on the constant
functions. This in turn leads to the effect that the associated
Laplacian is also not zero on the constant functions. For general
graphs without any geometric interpretation this is just a matter
of choice. However,  the choice of $d$ matters if one wants a
consistent continuum limit of the graph. One cannot expect
convergence of the graph Laplacian associated to the difference
operator $d$ of Equation \eqref{eq:denny-difference} towards a
Laplacian, since as each of the graph Laplacians in the sequence
is not zero on the constant functions, also the limit operator
will share this property unless $\lim_{n\rightarrow \infty}
d(X_i)=c, \forall i=1,\ldots,n$, where $c$ is a constant. We
derive also the limit operator of the graph Laplacian induced by
the difference operator of Equation \eqref{eq:denny-difference}
introduced by Zhou et al. in the machine learning literature and
usually denoted as the normalized graph Laplacian in spectral
graph theory \citep{Chung1997}.
\end{remark}
Obviously,  in the finite case $d$ is always a bounded operator.
The adjoint operator $d^*:\H(E,\phi) \rightarrow \H(V,\chi)$ is
defined by
\[ \inner{df,u}_E = \inner{f,d^*u}_V, \quad \forall \, f \in H(V,\chi), \quad u \in \H(E,\phi).\]
\begin{lemma}\label{le:difference-adjoint}
The adjoint $d^*:\H(E,\phi) \rightarrow \H(V,\chi)$ of the
difference operator $d$ is explicitly given by:
\begin{align} \label{eq:adjoint-dir}
(d^*u)(l) = \frac{1}{2\chi_l}\bigg(\frac{1}{n}\sum_{i=1}^n
\gamma(w_{il})\, u_{il}\, \phi(w_{il})
                                 - \frac{1}{n}\sum_{i=1}^n \gamma(w_{li})\, u_{li}\, \phi(w_{li})\bigg).
\end{align}
\end{lemma}
\begin{proof}
Using the indicator function $f(j)=\Id_{j=l}$ it is
straightforward to derive:
\begin{align*}
  \frac{1}{n}\chi_l \;(d^*u)(l)
&= \inner{d\Id_{\cdot=l},u}_E
= \frac{1}{2n^2}\sum_{i=1}^n \Big(\gamma(w_{il})u_{il}\phi(w_{il}) - 
  \gamma(w_{li})u_{li}\phi(w_{li})\Big),
\end{align*}
where we have used $\inner{d\Id_{\cdot=l},u}_E
 = \frac{1}{2n^2}\sum_{i,j=1}^n (d\Id_{\cdot=l})_{ij} u_{ij} \phi(w_{ij})$.
\end{proof}
The first term of the rhs of \eqref{eq:adjoint-dir} can be
interpreted as the outgoing flow, whereas the second term can be
seen as the ingoing flow. The corresponding continuous counterpart
of $d$ is the gradient of a function and for $d^*$ it is the
divergence of a vector-field, measuring the infinitesimal
difference between in- and outgoing flow.
\subsection{The General Graph Laplacian}\label{subsec:graph Laplacian}
\begin{definition}[graph Laplacian for a directed graph]
Given Hilbert spaces $\H(V,\chi)$\\ and $\H(E,\phi)$ and the
difference operator $d:\H(V,\chi) \rightarrow \H(E,\phi)$ the
graph Laplacian $\Delta: \H(V,\chi) \rightarrow \H(V,\chi)$ is
defined as
\[ \Delta =d^*d.\]
\end{definition}
\begin{lemma}\label{le:laplacian-explicit}
Explicitly, $\Delta: \H(V,\chi) \rightarrow \H(V,\chi)$ is given
as:
\begin{align}\label{def:laplacian-dir}
(\Delta f)(l)=\frac{1}{2 \chi_l}\bigg[ \frac{1}{n}\sum_{i=1}^n
\big(\gamma(w_{il})^2\phi(w_{il}) +
\gamma(w_{li})^2\phi(w_{li})\big)\big(f(l)-f(i)\big)
                 \bigg].
\end{align}
\end{lemma}
\begin{proof}
The explicit expression $\Delta$ can be easily derived by plugging
the expression of $d^*$ and $d$ together:
\begin{align*}
(d^*df)(l)=&\frac{1}{2 \chi_l }\bigg(\frac{1}{n}\sum_{i=1}^n
\gamma(w_{il})^2[f(l)-f(i)]\phi(w_{il}) -
            \frac{1}{n}\sum_{i=1}^n \gamma(w_{li})^2[f(i)-f(l)]\phi(w_{li})\bigg) \nonumber \\
          =&\frac{1}{2 \chi_l}\Big[ f(l)\frac{1}{n}\sum_{i=1}^n  \widehat{w}_{ij}
             -\frac{1}{n}\sum_{i=1}^n f(i)\widehat{w}_{ij}\Big],
\end{align*}
where we have introduced
$\widehat{w}_{ij}=\left(\gamma(w_{il})^2\phi(w_{il}) +
\gamma(w_{li})^2\phi(w_{li})\right)$.
\end{proof}
\begin{proposition}
$\Delta$ is self-adjoint and positive semi-definite.
\end{proposition}
\begin{proof}
By definition, $\inner{f,\Delta g}_V=\inner{df,dg}_E =
\inner{\Delta f,g}_V$, and $\inner{f,\Delta f}_V=\inner{df,df}_E
\geq 0.$
\end{proof}

\subsection{The Special Case of an Undirected Graph}\label{sec:undirected-graph}
In the case of an undirected graph we have $w_{ij}=w_{ji}$, that
is whenever there is an edge from $i$ to $j$ there is an edge with
the same value from $j$ to $i$. This implies that there is no
difference between in- and outgoing edges. Therefore, $d^{out}_i
\equiv d^{in}_i$, so that we will denote the degree function by
$d$ with $d_i=\frac{1}{n}\sum_{j=1}^n w_{ij}$. The same for the
weights in $H_V$, $\chi_{out}\equiv \chi_{in}$, so that we have
only one function $\chi$. If one likes to interpret functions on
$E$ as flows, it is reasonable to restrict the space $\H_E$ to
antisymmetric functions since symmetric functions are associated
to flows which transport the same mass from vertex $i$ to vertex
$j$ and back. Therefore, as a net effect, no mass is transported
at all so that from a physical point of view these functions
cannot be observed at all. Since anyway we consider only functions
on the edges of the form $df$ (where $f$ is in $\H_V$) which are
by construction antisymmetric, we will not do this restriction
explicitly. The adjoint $d^*$ simplifies in the undirected case to
\begin{equation} \label{eq:adjoint-undir}
 (d^*u)(l) = \frac{1}{2\chi(d_l)} \frac{1}{n}\sum_{i=1}^n \gamma(w_{il})\phi(w_{il})(u_{il}-u_{li}), \nonumber
\end{equation}
and the general graph Laplacian on an undirected graph has the
following form:
\begin{definition}[graph Laplacian for an undirected graph]
Given Hilbert spaces\\ $\H(V,\chi)$ and $\H(E,\phi)$ and the
difference operator $d:\H(V,\chi) \rightarrow \H(E,\phi)$ the
graph Laplacian $\Delta: \H(V,\chi) \rightarrow \H(V,\chi)$ is
defined as
\[ \Delta =d^*d. \]
Explicitly, for any vertex $l$, we have
\begin{align} \label{eq:laplacian-undir}
(\Delta f)(l)=(d^*df)(l)=\frac{1}{\chi(d_l)}\bigg[
f(l)\frac{1}{n}\sum_{i=1}^n  \gamma^2(w_{il})\phi(w_{il})
              -\frac{1}{n}\sum_{i=1}^n f(i)\gamma^2(w_{il})\phi(w_{il})\bigg].
\end{align}
\end{definition}
In the literature one finds the following special cases of the
general graph Laplacian. Unfortunately there exist no unique names
for the three graph Laplacians we introduce here, most of the time
all of them are just called graph Laplacians. Only the term
'unnormalized' or 'combinatorial' graph Laplacian seems to be
established now. However,  the other two could both be called
normalized graph Laplacian. Since the first one is closely related
to a random walk on the graph we call it random walk graph
Laplacian
and the other one normalized graph Laplacian.\\
The 'random walk' graph Laplacian is defined as:
\begin{align}\label{eq:laplacian-normalized}
(\Delrw f)(i)=f(i)-\frac{1}{d_i}\frac{1}{n}\sum_{j=1}^n
w_{ij}f(j), \quad \Delrw f = (\Id - D^{-1}W)f,
\end{align}
where the matrix $D$ is defined as $D_{ij}=d_i \, \delta_{ij}$.
Note that $P=D^{-1}W$ is a stochastic matrix and therefore can be
used to define a Markov random walk on $V$, see for example
\cite{Woe2000}. The 'unnormalized' (or 'combinatorial') graph
Laplacian is defined as
\begin{align}\label{eq:laplacian-unnormalized}
(\Delu f)(i)=d_i f(i) - \frac{1}{n}\sum_{j=1}^n w_{ij}f(j), \quad
\Delu f=(D-W)f.
\end{align}
We have the following conditions on $\chi,\gamma$ and $\phi$ in
order to get these Laplacians:
\begin{align*}
\forall \, e_{ij} \in E: \quad \mathrm{rw:} 
\quad
\frac{\gamma^2(w_{ij})\phi(w_{ij})}{\chi(d_i)}=\frac{w_{ij}}{d_i},
\quad
\mathrm{unnorm:} 
\quad \frac{\gamma^2(w_{ij})\phi(w_{ij})}{\chi(d_i)}=w_{ij}.
\end{align*}
We observe that by choosing $\Delrw$ or $\Delu$ the functions
$\phi$ and $\gamma$ are not fixed. Therefore it can cause
confusion if one speaks of the 'random walk' or 'unnormalized'
graph Laplacian without explicitly defining the
corresponding Hilbert spaces and the difference operator.\\
We also consider the normalized graph Laplacian $\Deln$ introduced
by \cite{Chung1997,Zhou2003} using the difference operator of
Equation \eqref{eq:denny-difference} and the general spaces
$H(V,\chi)$ and $\H(E,\phi)$. Following the scheme a
straightforward calculation shows the following:
\begin{lemma}
The graph Laplacian $\Delta_{\mathrm{norm}}=d^*d$ with the
difference operator $d$ from Equation \eqref{eq:denny-difference}
can be explicitly written as
\[ (\Deln f)(l)= \frac{1}{n\, \chi(d_l)\,\sqrt{d_l}}\bigg[ \frac{f(l)}{\sqrt{d_l}}\frac{1}{n}\sum_{i=1}^n  \gamma^2(w_{il})\phi(w_{il})
 -\frac{1}{n}\sum_{i=1}^n \frac{f(i)}{\sqrt{d_i}}\gamma^2(w_{il})\phi(w_{il})\bigg] \]
The choice $\chi(d_l)=1$ and $\gamma^2(w_{il})\phi(w_{il})=w_{il}$
leads then to the graph Laplacian proposed in
\cite{ChuLan1996,Zhou2003},
\[  (\Deln f)(l)= \frac{1}{n\,\,\sqrt{d_l}}\bigg[ \frac{f(l)}{\sqrt{d_l}}d_l
 -\frac{1}{n}\sum_{i=1}^n \frac{f(i)}{\sqrt{d_i}}w_{li}\bigg]=\frac{1}{n}\bigg[ f(l)
 -\frac{1}{n}\sum_{i=1}^n f(i)\frac{w_{il}}{\sqrt{d_l\,d_i}}\bigg], \]
or equivalently
\[ \Deln f = D^{-\frac{1}{2}} (D-W) D^{-\frac{1}{2}}f=(\Id - D^{-\frac{1}{2}}W D^{-\frac{1}{2}})f. \]
\end{lemma}
Note that $\Delu=D \Delrw$ and $\Deln=D^{-\frac{1}{2}}\Delu
D^{-\frac{1}{2}}$.


\section{Limit of the Graph Laplacian for Random Neighborhood Graphs}\label{sec:limit}
Before we state the convergence results for the three graph
Laplacians on random neighborhood graphs, we first have to define
the limit operator. Maybe not surprisingly, in general the
Laplace-Beltrami operator will \emph{not} be the limit operator of
the graph Laplacian. Instead it will converge to the weighted
Laplace-Beltrami operator which is the natural generalization of
the Laplace-Beltrami operator for a Riemannian manifold equipped
with a non-uniform probability measure. The definition of this
limit operator and a discussion of its use for different
applications in clustering, semi-supervised learning and
regression is the topic of the next section, followed by a sketch
of the convergence results.

\subsection{Construction of the Neighborhood Graph}\label{sec:construction}
We assume to have a sample $X_i, \, i=1,\ldots,n$ drawn i.i.d.
from a probability measure $P$ which has support on a submanifold
$M$. For the exact assumptions regarding $P$, $M$ and the kernel
function $k$ used to define the weights we refer to Section
\ref{sec:assumptions}. The sample then determines the set of
vertices $V$ of the graph. Additionally we are given a certain
kernel function $k:\R_+ \rightarrow \R_+$ and the neighborhood
parameter $h \in \R^*_+$. As proposed by
\cite{Lafon2004,CoiLaf2005}, we use this kernel function $k$ to
define the following family of data-dependent kernel functions
$\tilde{k}_{\lambda,h}$ parameterized by $\lambda \in \R$ as:
\[ \tilde{k}_{\lambda,h}(X_i,X_j)= \frac{1}{h^m}\frac{k(\norm{X_i-X_j}^2/h^2)}{[d_{h,n}(X_i)d_{h,n}(X_j)]^\lambda}, \]
where $d_{h,n}(X_i)=\frac{1}{n}\sum_{i=1}^n \frac{1}{h^m}
k(\norm{X_i-X_j}^2/h^2)$ is the degree function introduced in
Section \ref{sec:graphs} with respect to the edge-weights
$\frac{1}{h^m} k(\norm{X_i-X_j}^2/h^2)$. Finally we use
$\tilde{k}_{\lambda,h}$ to define the weight $w_{ij}=w(X_i,X_j)$
of the edge between the points $X_i$ and $X_j$ as
\[ w_{\lambda,h}(X_i,X_j) = \tilde{k}_{\lambda,h}(X_i,X_j).\]
Note that the case $\lambda=0$ corresponds to weights with no
data-dependent modification. The parameter $h \in \R^*_+$
determines the neighborhood of a point since we will assume that
$k$ has compact support, that is $X_i$ and $X_j$ have an edge if
$\norm{X_i-X_j}\leq hR_k$ where $R_k$ is the support of kernel
function. Note that we will have $k(0)=0$, so that there are no
loops in the graph. This assumption is not
necessary, but it simplifies the proofs and makes some of the estimators unbiased.\\
In Section \ref{sec:undirected-graph} we introduced the random
walk, the unnormalized and the normalized graph Laplacian. From
now on we consider these graph Laplacians for the random
neighborhood graph, that is the weights of the graph $w_{ij}$ have
the form $w_{ij}= w(X_i,X_j)=\tk_{\lambda,h}(X_i,X_j)$. Using the
kernel function we can easily extend the graph Laplacians to the
whole submanifold $M$. These extensions can be seen as estimators
for the Laplacian on $M$. We introduce also the extended degree
function $\td_{\lambda,h,n}$ and the average operator
$\tA_{\lambda,h,n}$,
\begin{align*}
    \td_{\lambda,h,n}(x)=\frac{1}{n}\sum_{j=1}^n \tilde{k}_{\lambda,h}(x,X_j), \quad \quad
 (\tA_{\lambda,h,n}f)(x)=\frac{1}{n}\sum_{j=1}^n \tilde{k}_{\lambda,h}(x,X_j)f(X_j).
\end{align*}
Note that $\td_{\lambda,h,n}=\tA_{\lambda,h,n}1$. The extended
graph Laplacians are defined as follows:
\begin{align}\label{eq:graph-laplacian-final}
\hspace{-0.5cm}\mathrm{random\, walk} &&(\Delrw_{\lambda,h,n}f)(x)
&=\frac{1}{h^2}\Big(f- \tfrac{1}{\td_{\lambda,h,n}}
\tA_{\lambda,h,n}f\Big)(x)
                           =\frac{1}{h^2}\left(\frac{\tA_{\lambda,h,n}g}{\td_{\lambda,h,n}} \right)(x),  \\
\hspace{-0.5cm}\mathrm{unnormalized} &&(\Delu_{\lambda,h,n}f)(x) &= \frac{1}{h^2}\Big(\td_{\lambda,h,n}f - \tA_{\lambda,h,n}f\Big)(x) 
                            = \frac{1}{h^2} (\tA_{\lambda,h,n}g)(x),                                 \\
\hspace{-0.5cm}\mathrm{normalized}&&(\Deln_{\lambda,h,n}f)(x) &= \tfrac{1}{h^2\,\sqrt{\td_{\lambda,h,n}(x)}}\bigg(\td_{\lambda,h,n}\tfrac{f}{\sqrt{\td_{\lambda,h,n}}} - \Big(\tA_{\lambda,h,n}\tfrac{f}{\sqrt{\td_{\lambda,h,n}}}\Big)\bigg)(x) \nonumber \\
\hspace{-0.5cm}&&         &=
\tfrac{1}{h^2\,\sqrt{\td_{\lambda,h,n}(x)}}(\tA_{\lambda,h,n}g')(x)
,
\end{align}
where we have introduced $g(y):=f(x)-f(y)$ and
$g'(y):=\frac{f(x)}{\sqrt{\td_{\lambda,h,n}(x)}}-\frac{f(y)}{\sqrt{\td_{\lambda,h,n}(y)}}$.
Note that all extensions reproduce the graph Laplacian on the
sample:
\[  (\Delta f)(i)=(\Delta f)(X_i)=(\Delta_{\lambda,h,n}f)(X_i), \quad \forall i=1,\ldots,n.\]
The factor $1/h^2$ arises by introducing a factor $1/h$ in the
weight $\gamma$ of the derivative operator $d$ of the graph. This
is necessary since $d$ is supposed to approximate a derivative.
Since
the Laplacian corresponds to a second derivative we get from the definition of the graph Laplacian a factor $1/h^2$.\\
We would like to note that in the case of the random walk and and
the normalized graph Laplacian the normalization with $1/h^m$ in
the weights cancels out, whereas it does not cancel for the
unnormalized graph Laplacian  except in the case $\lambda=1/2$.
The problem here is that in general the intrinsic dimension $m$ of
the manifold is unknown. Therefore a normalization with the
correct factor $\frac{1}{h^m}$ is not possible, and in the limit
$h \rightarrow 0$ the estimate of the unnormalized graph Laplacian
will generally either vanish or blow up. The easy way to
circumvent this is just to rescale the whole estimate such that
$\frac{1}{n}\sum_{i=1}^n \td_{\lambda,h,n}(X_i)$ equals a fixed
constant for every $n$. The disadvantage is that this method of
rescaling introduces a global factor in the limit. A more elegant
way might be to  simultaneously estimate the dimension $m$ of the
submanifold and use the estimated dimension to calculate the
correct normalization factor, see e.g. \cite{HeiAud2005}. However,
in this work we assume for simplicity that for the unnormalized
graph Laplacian the intrinsic dimension $m$ of the submanifold is
known. It might be interesting to consider both estimates
simultaneously, but we leave this as an open problem.\\
We will consider in the following the limit $h\rightarrow 0$, that
is the neighborhood of each point $X_i$ shrinks to zero. However,
since $n \rightarrow \infty$ and $h$ as a function of $n$
approaches zero sufficiently slow, the number of points in each
neighborhood approaches $\infty$, so that roughly spoken sums
approximate the corresponding integrals. This is the basic
principle behind our convergence result and is well known in the
framework of nonparametric regression \citep[see][]{Gyo04}.

\subsection{The Weighted Laplacian and the Continuous Smoothness Functional}
The Laplacian is one of the most prominent operators in
mathematics. The following general properties are taken from the
books of \cite{Rosenberg1997} and 
\cite{Berard1986}. It occurs in many partial differential
equations governing physics, mainly because in Euclidean space it
is  the only linear second-order differential operator which is
translation and rotation invariant. In Euclidean space $\R^d$ it
is defined as $  \Delta_{\R^d}f=\div(\grad f)=\sum_{i=1}^d
\partial^2_i f .$ Moreover,  for any domain $\Omega \subseteq
\R^d$ it is a negative-semidefinite symmetric operator on
$C^\infty_c(\Omega)$, which is a dense subset of $L_2(\Omega)$
(formally self-adjoint), and satisfies
\[ \int_\Omega f \Delta h \,dx = - \int_\Omega \inner{\nabla f,\nabla h} dx.\]
It can be extended to a self-adjoint operator on $L_2(\Omega)$ in
several ways depending on the choice of boundary conditions. For
any compact domain $\Omega$ (with suitable boundary conditions) it
can be shown that $\Delta$ has a pure point spectrum and the
eigenfunctions are smooth and form
a complete orthonormal basis of $L_2(\Omega)$, see e.g. \cite{Berard1986}.\\
The Laplace-Beltrami operator on a Riemannian manifold $M$ is the
natural equivalent of the Laplacian in $\R^d$, defined as
\[ \Delta_M f = \div(\grad f) = \nabla^a \nabla_a f. \]
However,  the more natural definition is the following. For any
$f,g \in C^\infty_c(M)$, we have
\[ \int_M f \Delta h\, dV(x) = - \int_M \inner{\nabla f,\nabla h} dV(x), \]
where $dV=\sqrt{\det g}\, dx$ is the natural volume element of
$M$. This definition allows easily an extension to the case where
we have a Riemannian manifold $M$ with a measure $P$ . In this
paper $P$ will be the probability measure generating the data. We
assume in the following that $P$ is absolutely continuous wrt the
natural volume element $dV$ of the manifold. Its density is
denoted by $p$. Note that the case when the probability measure is
absolutely continuous wrt the Lebesgue measure on $\R^d$ is a
special case of our setting.\\
A recent review article about the weighted Laplace-Beltrami
operator is \citep{Gri2006}.
\begin{definition}[Weighted Laplacian]
Let $(M,g_{ab})$ be a Riemannian manifold\ with measure $P$ where
$P$ has a differentiable and positive density $p$ with respect to
the natural volume element $dV=\sqrt{\det g}\,dx$, and let
$\Delta_M$ be the Laplace-Beltrami operator on $M$. For $s\in\R$,
we define the $s$-th weighted Laplacian $\Delta_{s}$ as
\begin{align}\label{def:weighted-Laplacian}
\Delta_{s} := \Delta_M + \frac{s}{p} g^{ab}(\nabla_a p) \nabla_b =
\frac{1}{p^s}g^{ab}\nabla_a (p^s
\nabla_b)=\frac{1}{p^s}\div(p^s\,\grad).
\end{align}
\end{definition}
This definition is motivated by the following equality, for $f,g
\in C^\infty_c(M)$,
\begin{align}\label{eq:laplacian-self-adjoint}
\int_M f (\Delta_s g)\, p^s dV = \int_M f \big(\Delta g +
\frac{s}{p}\inner{\nabla p,\nabla g}\big) p^s dV=-\int_M
\inner{\nabla f,\nabla g} p^s dV,
\end{align}
The family of weighted Laplacians contains two cases which are
particularly interesting. The first one, $s=0$, corresponds to the
standard Laplace-Beltrami operator. This case is interesting if
one only wants to use properties of the geometry of the manifold
but not of the data generating probability measure. The second
case, $s=1$, corresponds to the standard
weighted Laplacian $\Delta_1 = \frac{1}{p}\nabla^a(p \nabla_a)$.\\
In the next section it will turn out that through a data-dependent
change of the weights of the graph we can get the just defined
weighted Laplacians as the limit operators of the graph Laplacian.
The rest of this section will be used to motivate the importance
of the understanding of this limit in different applications.
Three different but closely related properties of the Laplacian
are used in machine learning
\begin{itemize}
\item The Laplacian generates the diffusion process. In
semi-supervised learning algorithms with a small number of labeled
points one would like to propagate the labels along regions of
high-density, see \cite{ZhuGha2002, ZhuGhaLaf2003}. The limit
operator $\Delta_s$ shows the influence of a non-uniform density
$p$. The second term $ \frac{s}{p}\inner{\nabla p,\nabla f}$ leads
to an anisotropy in the diffusion. If $s<0$ this term enforces
diffusion in the direction of the maximum of the density whereas
diffusion in the direction away from the maximum of the density is
weakened. If $s<0$ this is just the other way round. \item The
smoothness functional induced by the weighted Laplacian
$\Delta_s$, see Equation \ref{eq:laplacian-self-adjoint}, is given
by
\[ S(f)=\int_M \norm{\nabla f}^2 p^s \,dV. \]
For $s>0$ this smoothness functional prefers  functions which are
smooth in high-density regions whereas unsmooth behavior in
low-density is penalized less. This property can also be
interesting in semi-supervised learning where one assumes
especially when only a few labeled points are known that the
classifier should be constant in high-density regions whereas
changes of the classifier are allowed in low-density regions, see
\cite{BouChaHei2003} for some discussion of this point and
\cite{Hei2005,Hei2006} for a proof of convergence of the
regularizer induced by the graph Laplacian towards the smoothness
functional $S(f)$. In Figure \ref{fig:density-manifold} this is
illustrated by mapping a density profile in $\R^2$ onto a
two-dimensional manifold.
\begin{figure}[h]
$$\epsfig{file=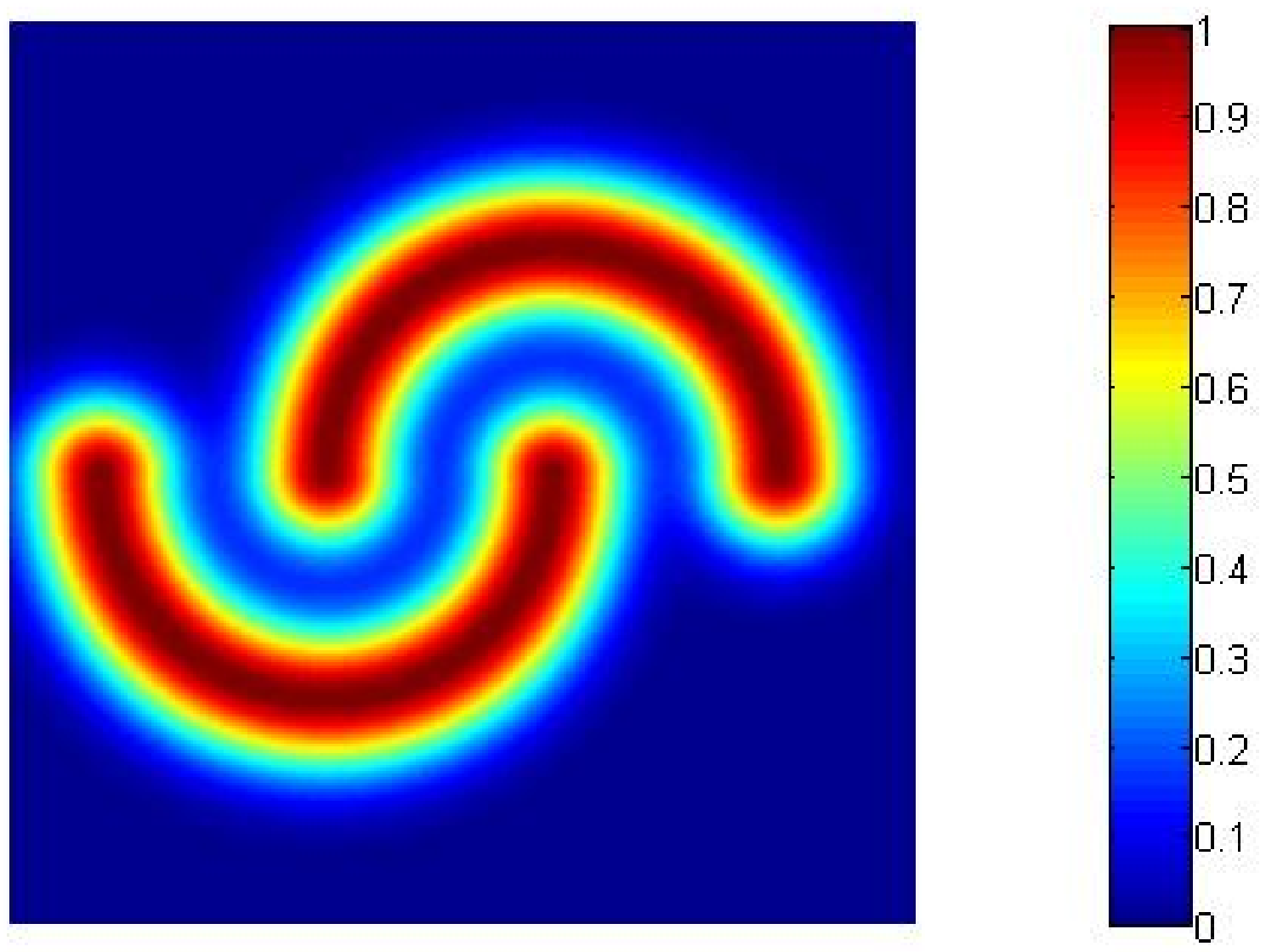,width=6cm} \quad \quad  \epsfig{file=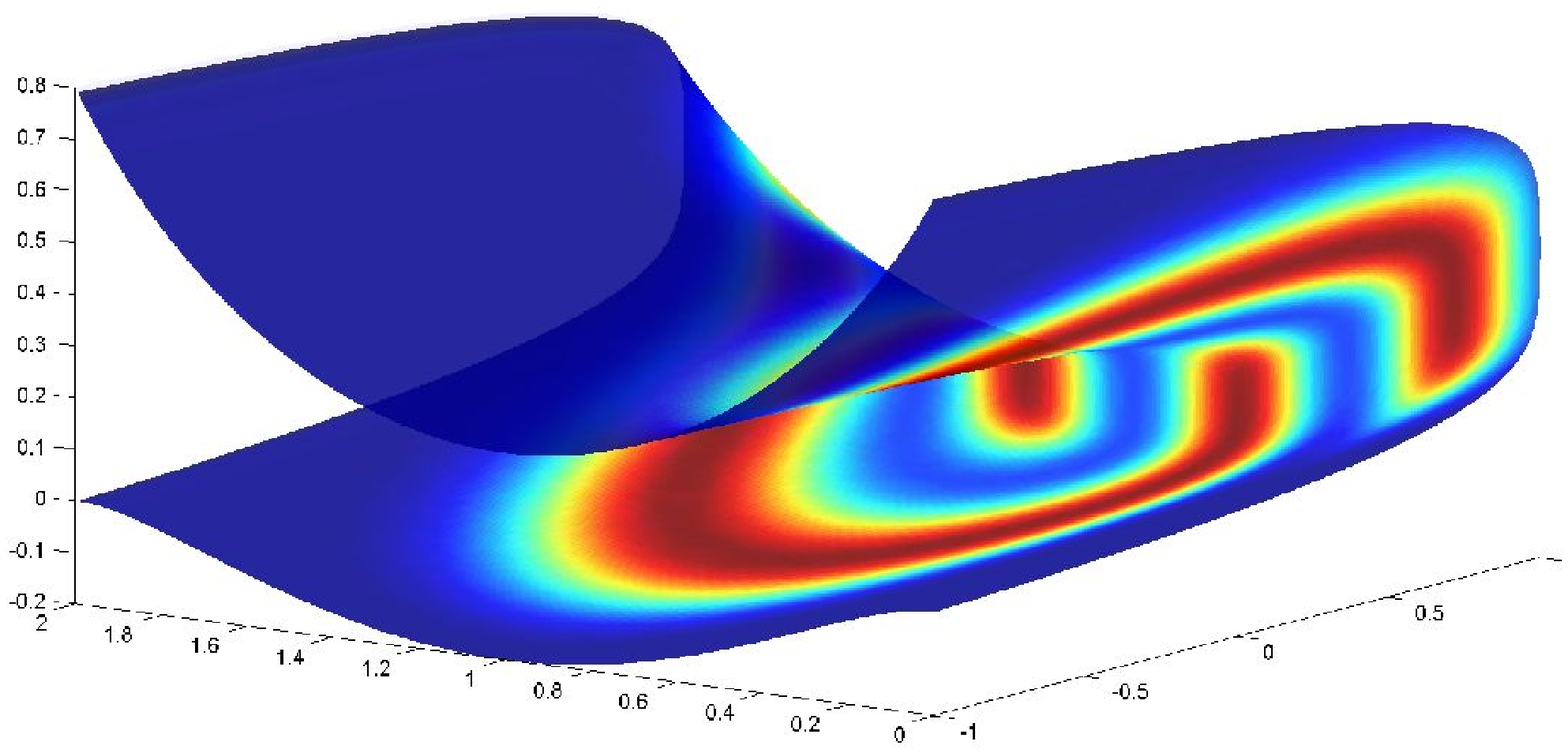,width=7cm}$$
\caption{A density profile mapped onto a two-dim. submanifold in
$\R^3$ with two clusters.} \label{fig:density-manifold}
\end{figure}
However,  also the case $s<0$ can be interesting. Minimizing the
smoothness functional $S(f)$ implies that one enforces smoothness
of the function $f$ where one has little data, and one allows the
function to vary more where one has sampled a lot of data points.
Such a penalization has been considered by \cite{CanEli1999} for
regression.

\item The eigenfunctions of the Laplacian $\Delta_s$ can be seen
as the limit partioning of spectral clustering for the normalized
graph Laplacian (however,  a rigorous mathematical proof has not
been given yet, see \cite{LuxBelBou2004} for a convergence result
for fixed $h$). If $s=0$ one gets just a geometric clustering in
the sense that irrespectively of the probability measure
generating the data the clustering is determined by the geometry
of the submanifold. If $s>0$ the eigenfunction corresponding to
the first non-zero eigenvalue is likely to change its sign in a
low-density region. This argument follows from the previous
discussion on the smoothness functional $S(f)$ and the
Rayleigh-Ritz principle. Let us assume for a moment that $M$ is
compact without boundary and that $p(x)>0, \forall \, x \in M$,
then the eigenspace corresponding to the first eigenvalue
$\lambda_0=0$ is given by the constant functions. The first
non-zero eigenvalue can then be determined by the Rayleigh-Ritz
variational principle
\[ \lambda_1 = \inf_{ u \in C^\infty(M)} \bigg\{ \frac{\int_M \norm{\nabla u}^2 p(x)^s dV(x)}{\int_M u^2(x) p(x)^s dV(x)} \; \Big| \;
  \int_M u(x) \, p(x)^s dV(x) = 0 \bigg \}. \]
Since the first eigenfunction has to be orthogonal to the constant
functions, it has to change its sign. However,  since
$\norm{\nabla u}^2$ is weighted by a power of the density $p^s$ it
is obvious that for $s>0$ the function will change its sign in a
region of low density.
\end{itemize}

\subsection{Limit of the Graph Laplacians}\label{sec:limit-sketch}
The following theorem summarizes and slightly weakens the results
of Theorem \ref{th:limit-laplacian} and Theorem
\ref{th:limit-laplacian-normalized} of Section
\ref{sec:continuum-limit}.\vspace{+2mm}\\
{\textbf{Main Result \;}{\it Let $M$ be a $m$-dimensional
submanifold in $\R^d$, $\{X_i\}_{i=1}^n$ a sample from a
probability measure $P$ on $M$ with density $p$. Let $x \in
M\backslash \partial M$ and define $s=2(1-\lambda)$. Then under
technical conditions on the submanifold $M$, the kernel $k$ and
the density $p$ introduced in Section \ref{sec:continuum-limit},
if $h \rightarrow 0$ and $nh^{m+2} / \log n \rightarrow \infty$,}
\begin{align*}
\textrm{random walk:}             &&\lim_{n \rightarrow \infty}
\,(\Delrw_{\lambda,h,n}f)(x)
&\;\sim\; -(\Delta_{s} f)(x) &&\textrm{almost surely},\\
\textrm{unnormalized:}   &&\lim_{n \rightarrow \infty}
\,(\Delu_{\lambda,h,n}f)(x) &\;\sim\;
-p(x)^{1-2\lambda}\,(\Delta_{s} f)(x) &&\textrm{almost surely}.
\end{align*}
The optimal rate is obtained for $h(n)=O\Big( (\log n /
n)^\frac{1}{m+4}\Big)$. If $h \rightarrow 0$ and $nh^{m+4} / \log
n \rightarrow \infty$,
\begin{align*}
\textrm{normalized:}     &&\lim_{n \rightarrow \infty}\,
(\Deln_{\lambda,h,n}f)(x) &\;\sim\; -p(x)^{\frac{1}{2}-\lambda}
\Delta_{s} \bigg(\frac{f}{p^{\frac{1}{2}-\lambda}}\bigg)(x)
&&\textrm{almost surely}.
\end{align*}
{\it where $\sim$ means that there exists a constant only
depending on the kernel $k$ and $\lambda$ such that equality
holds.\vspace{+2mm}\\}} The first observation is that the
conjecture that the graph Laplacian approximates the
Laplace-Beltrami operator is only true for the uniform measure,
where $p$ is constant. In this case all limits agree up to
constants. However,  big differences arise when one has a
non-uniform measure on the submanifold, which is the generic case
in machine learning applications. In this case all limits disagree
and only the random walk graph Laplacian converges towards the
weighted Laplace-Beltrami operator which is the natural
generalization of the Laplace-Beltrami operator when the manifold
is equipped with a non-uniform probability measure. The
unnormalized graph Laplacian has the additional factor
$p^{1-2\lambda}$. However,  this limit is actually quite useful,
when one thinks of applications of so called label propagation
algorithms in semi-supervised learning. If one uses this graph
Laplacian as the diffusion operator to propagate the labeled data,
it means that the diffusion for $\lambda<1/2$ is faster in regions
where the density is high. The consequence is that labels in
regions of high density are propagated faster than labels in
low-density regions. This makes sense since  under the cluster
assumption labels in regions of low density are less informative
than labels in regions of high density. In general, from the
viewpoint of a diffusion process the weighted Laplace-Beltrami
operator $\Delta_s=\Delta_M + \frac{s}{p}\nabla p \nabla$ can be
seen as inducing an anisotropic diffusion due to the extra term
$\frac{s}{p}\nabla p \nabla $, which is directed towards or away
from increasing density depending on $s$. This is a desired
property in semi-supervised learning, where one actually wants
that the diffusion is mainly along
regions of the same density level in order to fulfill the cluster assumption.\\
The second observation is that the data-dependent modification of
edge weights allows to control the influence of the density on the
limit operator as observed by \cite{CoiLaf2005}. In fact one can
even eliminate it for $s=0$ resp. $\lambda=1$ in the case of the
random walk graph Laplacian. This could be interesting in computer
graphics where the random walk graph Laplacian is used for mesh
and point cloud processing, see e.g. \cite{Sor2005}. If one has
gathered points of a curved object with a laser scanner it is
likely that the points have a non-uniform distribution on the
object. Its surface is a two-dimensional submanifold in $\R^3$. In
computer graphics the non-uniform measure is only an artefact of
the sampling procedure and one is only interested in the
Laplace-Beltrami operator to infer geometric properties. Therefore
the elimination of the influence of a non-uniform measure on the
submanifold is of high interest there. We note that up to a
multiplication with the inverse of the density the elimination of
density effects is also possible for the unnormalized graph
Laplacian, but not for the normalized graph Laplacian. All
previous observations are naturally also true if the data does not
lie on a submanifold but has $d$-dimensional support
in $\R^d$.\\
The interpretation of the limit of the normalized graph Laplacian
is more involved. An expansion of the limit operator shows the
complex dependency on the density $p$:
\begin{align*}
p^{\frac{1}{2}-\lambda} \Delta_{s}
\bigg(\frac{f}{p^{\frac{1}{2}-\lambda}}\bigg)
&=\Delta_M f + \frac{1}{p}\nabla p \nabla f -
(\lambda-\frac{1}{2})^2 \,\frac{f}{p}\norm{\nabla p}^2 +
(\lambda-\frac{1}{2})\frac{f}{p}\Delta_M p
\end{align*}
We leave it to the reader to think of possible applications of this Laplacian.\\
The discussion shows that the choice of the graph Laplacian
depends on what kind of problem one wants to solve. Therefore, in
our opinion there is no universal best choice between the random
walk and the unnormalized graph Laplacian from a machine learning
point of view. However,  from a mathematical point of view only
the random walk graph Laplacian has the correct (pointwise) limit
to the weighted Laplace-Beltrami operator.


\section{Illustration of the Results} \label{sec:illustration}
In this section we want to illustrate the differences between the
three graph Laplacians and the control of the influence of the
data-generating measure via the parameter $\lambda$.
\subsection{Flat Space $\R^2$}
\begin{figure}
$$\epsfig{file=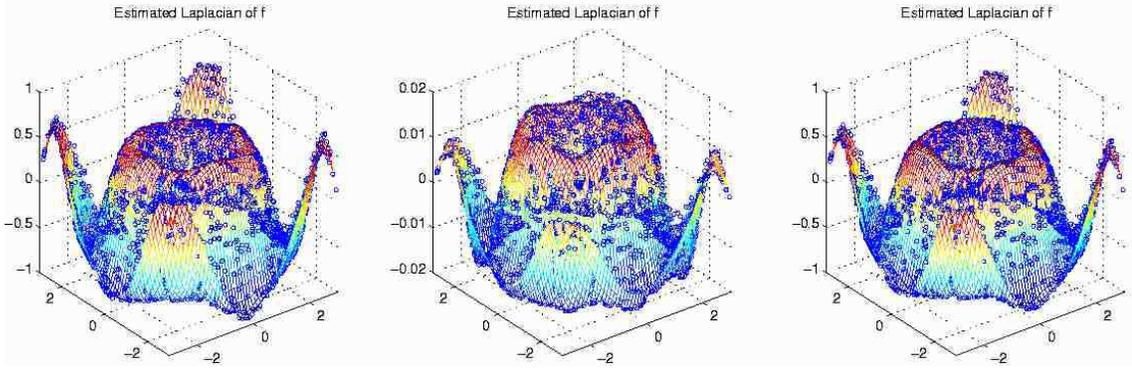,width=15cm}$$
\caption{For the uniform distribution all graph Laplacians,
$\Delrw_{\lambda,h,n}$, $\Delu_{\lambda,h,n}$ and
$\Deln_{\lambda,h,n}$ (from left to right) agree up to constants
for all $\lambda$. In the figure the estimates of the Laplacian
are shown for the uniform measure on $[-3,3]^2$ and the function
$f(x)=\sin(\frac{1}{2}\norm{x}^2)/\norm{x}^2$ with 2500 samples
and $h=1.4$.
\label{fig:illustration-uniform}}
\end{figure}

In the first example the data lies in Euclidean space $\R^2$. Here
we want to show two things. First, the sketch of the main result
shows that if the data generating measure is uniform all graph
Laplacians converge for all values of the reweighting parameter
$\lambda$ up to constants to the Laplace-Beltrami operator, which
is in the case of $\R^2$ just the standard Laplacian. In Figure
\ref{fig:illustration-uniform} the estimates of the three graph
Laplacians are shown for the uniform measure $[-3,3]^2$ and
$\lambda=0$. It can be seen that up to a scaling all estimates
agree very well.
In a second example we study the effect of a non-uniform
data-generating measure. In general all estimates disagree in this
case. We illustrate this effect in the case of $\R^2$ with a
Gaussian distribution $\mathcal{N}(0,1)$ as data-generating
measure and the simple function $f(x)=\sum_i x_i - 4$. Note that
$\Delta f=0$ so that for the random walk and the unnormalized
graph Laplacian only the anisotropic part of the limit operator,
$\frac{1}{p}\nabla p \nabla f$ is non-zero. Explicitly the limits
are given as
\begin{align*}
\Delrw &\sim& \Delta_s f &=\Delta f + \frac{s}{p} \nabla p \nabla f = - s \sum\nolimits_i x_i,\\
\Delu  &\sim& p^{1-2\lambda}\Delta_s f & = -s\, e^{-\frac{1-2\lambda}{2}\norm{x}^2}\sum\nolimits_i x_i,\\
\Deln  &\sim& p^{\frac{1}{2}-\lambda} \Delta_s
\tfrac{f}{p^{\frac{1}{2}-\lambda}}&
=  - \sum\nolimits_i x_i - \big(\sum\nolimits_i x_i -4\big)\Big[
(\lambda-\frac{1}{2}) (\frac{3}{2}-\lambda)\norm{x}^2
-2(\lambda-\frac{1}{2})\Big]
\end{align*}
This shows that even applied to simple functions there can be
large differences between the different limit operators provided
the samples come from a non-uniform probability measure. Note that
like in nonparametric kernel regression the estimate is quite bad
at the boundary. This well known boundary effect arises since at
the boundary one does not average over a full ball but only over
some part of a ball. Thus the first derivative $\nabla f$ of order
$O(h)$ does not cancel out so that multiplied with the factor
$1/h^2$ we have a term of order $O(1/h)$ which blows up. Roughly
spoken this effect takes place at all points of order $O(h)$ away
from the boundary, see also \citep{CoiLaf2005}.
\begin{center}
\begin{figure}
$$\epsfig{file=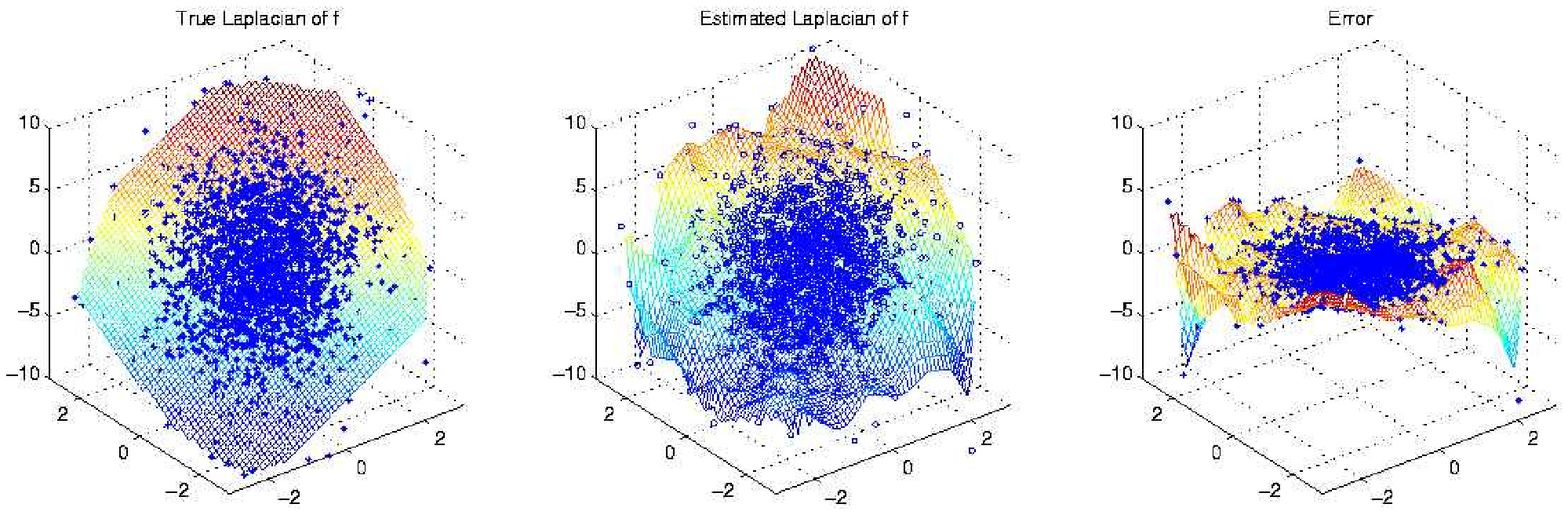,width=15cm}$$
$$\epsfig{file=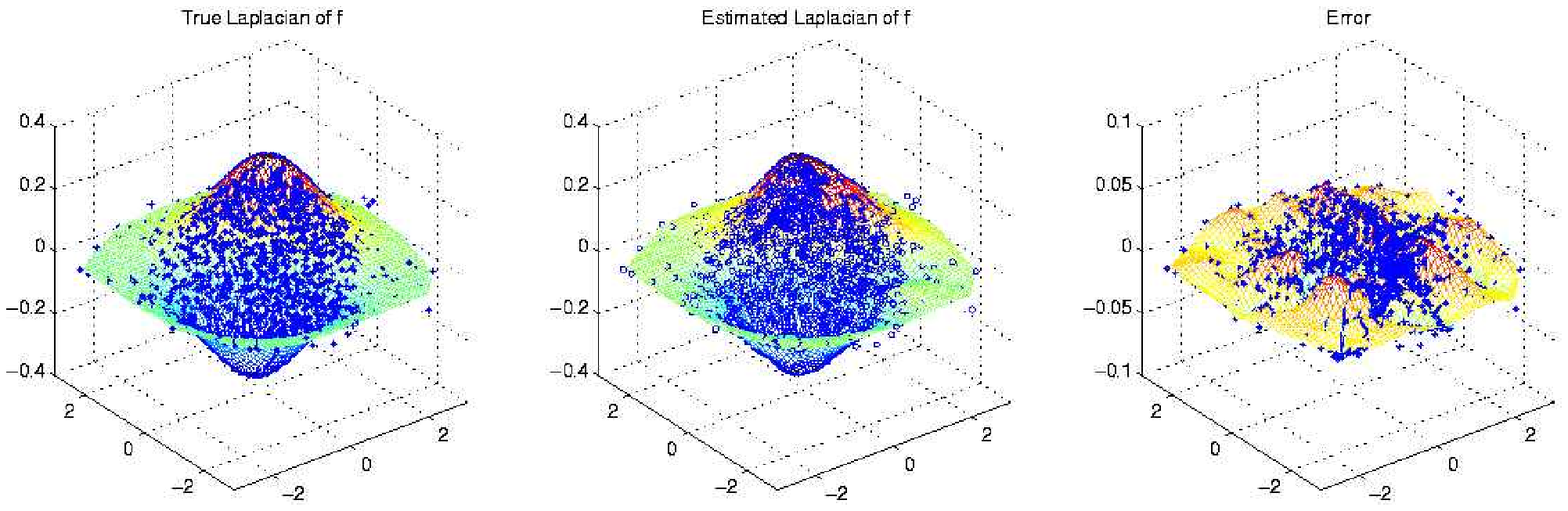,width=15cm}$$
$$\epsfig{file=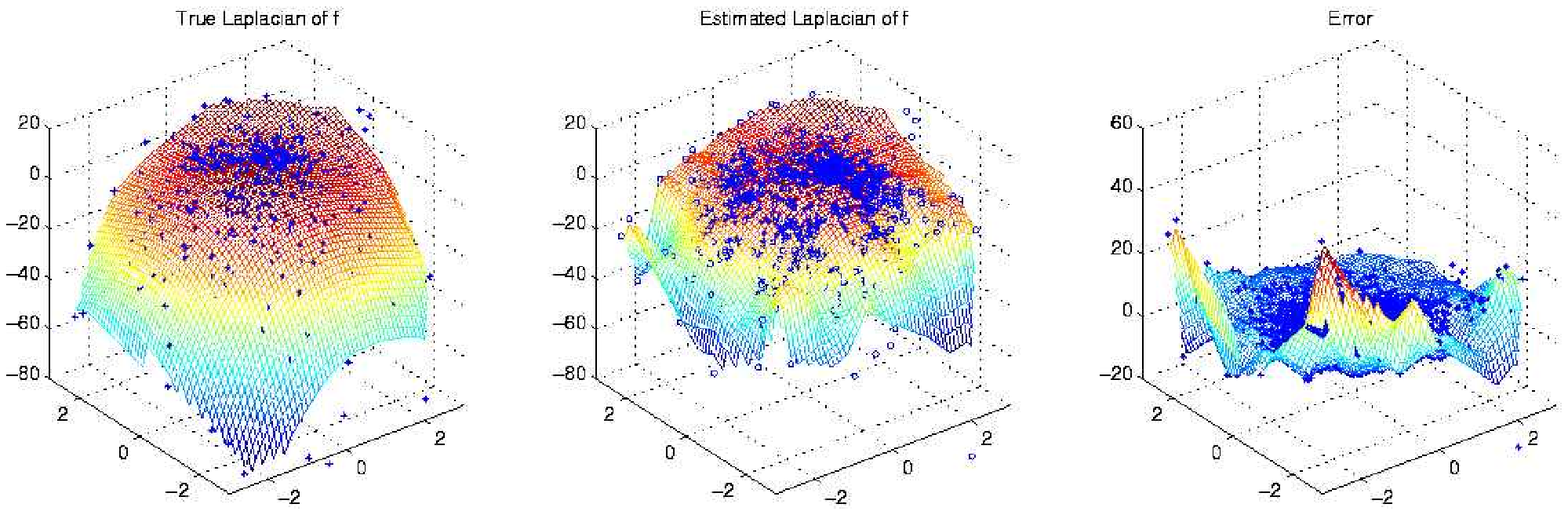,width=15cm}$$
\caption{Illustration of the differences of the three graph
Laplacians, random walk, unnormalized and normalized (from the
top) for $\lambda=0$. The function $f$ is $f=\sum_{i=1}^2 x_i -4$
and the 2500 samples come from a standard Gaussian distribution on
$\R^2$. The neighborhood size $h$ is set to $1.2$.
\label{fig:illustration-gaussian}}
\end{figure}
\end{center}

\subsection{The Sphere $S^2$}
\begin{center}
\begin{figure}
$$\epsfig{file=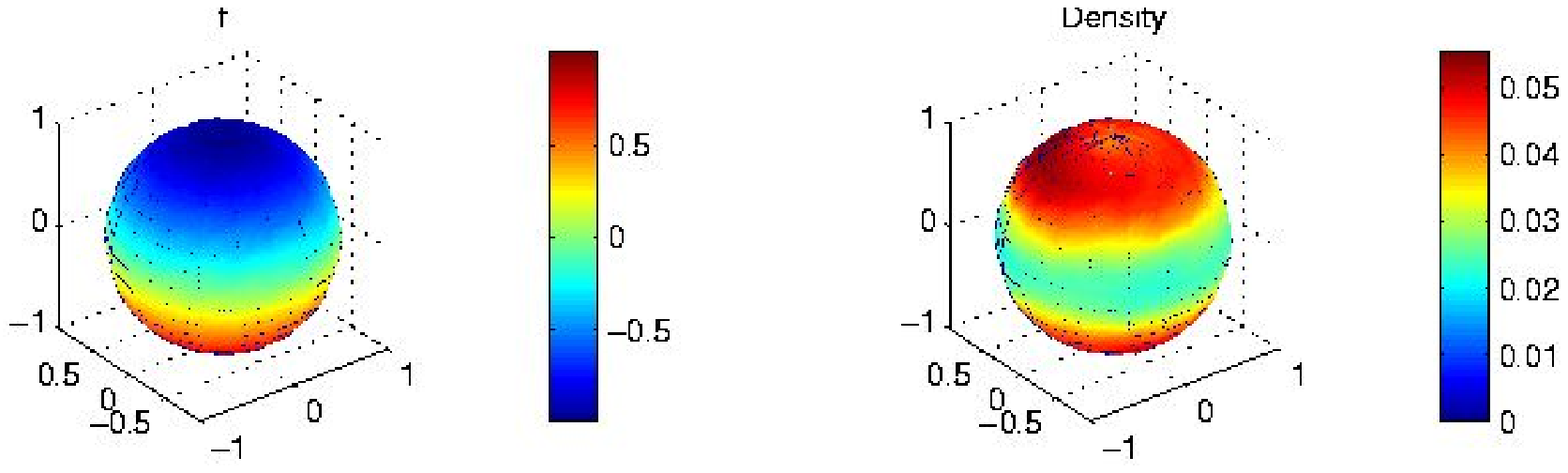,width=15cm}$$
$$\epsfig{file=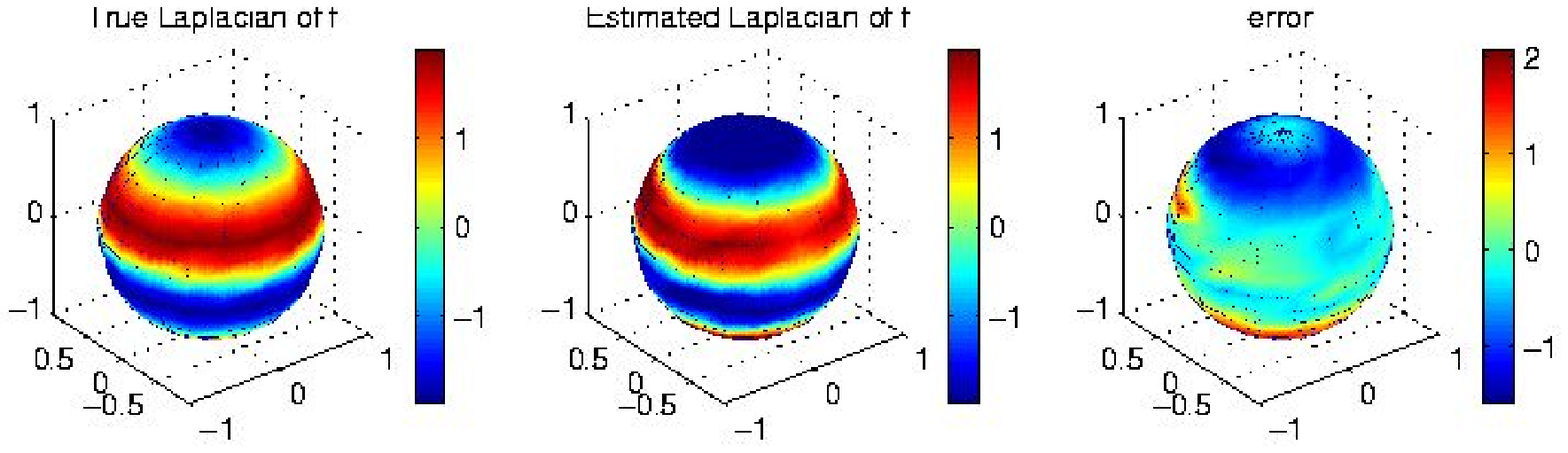,width=15cm}$$
$$\epsfig{file=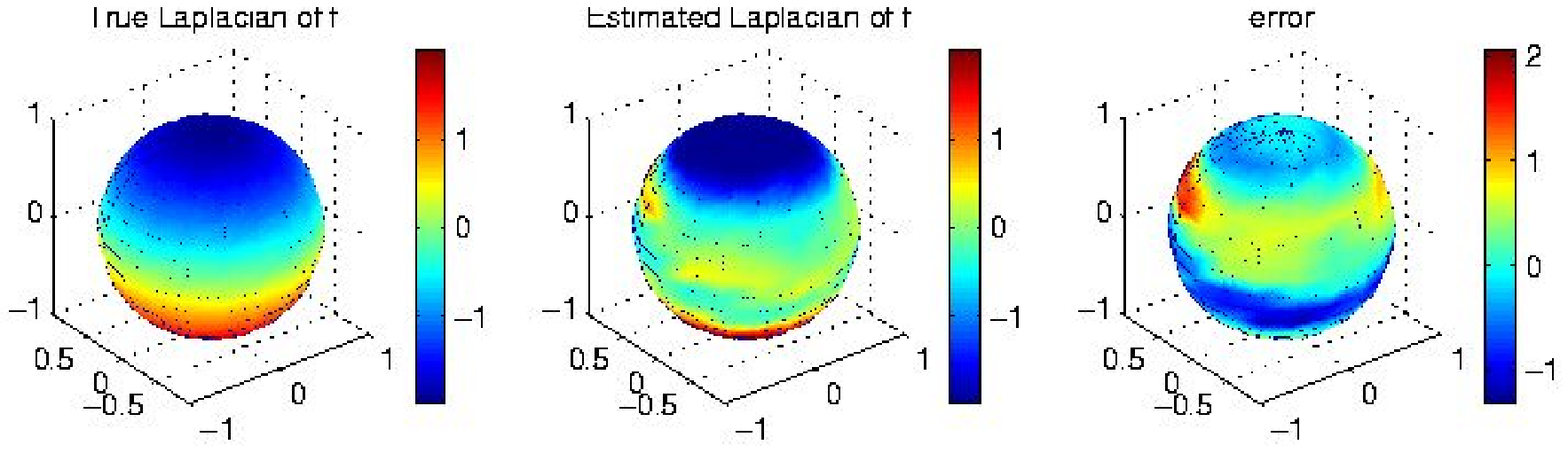,width=15cm}$$
$$\epsfig{file=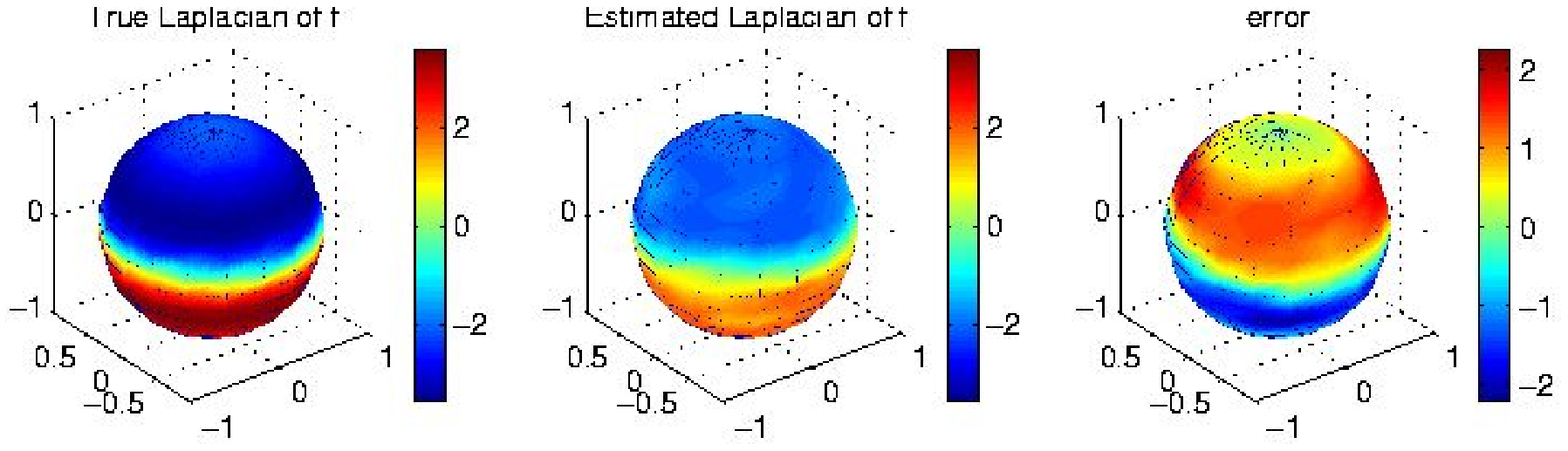,width=15cm}$$
\caption{Illustration of the effect of $\lambda=0,1,2$ (row $2-4$)
resulting in $s=-2,0,2$ for the sphere with a non-uniform
data-generating probability measure and the function
$f(\theta,\phi)=\cos(\theta)$ (row $1$) for the random walk
Laplacian
with $n=2500$ and $h=0.6$
\label{fig:illustration-sphere}}
\end{figure}
\end{center}

In our next example we consider the case where the data lies on a
submanifold $M$ in $\R^d$. Here we want to illustrate in the case
of a sphere $S^2$ in $\R^3$ the control of the influence of the
density via the parameter $\lambda$.  In this case we sample from
the probability measure with density
$p(\phi,\theta)=\frac{1}{8\pi}+\frac{3}{8\pi}\cos^2(\theta)$ in
spherical coordinates with respect to the volume element
$dV=\sin(\theta)d\theta d\phi$. This density has a two-cluster
structure on the sphere, where the northern and southern
hemisphere represent one cluster. An estimate of the density $p$
is shown in the Figure \ref{fig:illustration-sphere}. We show the
results of the random walk graph Laplacian together with the
result of the weighted Laplace-Beltrami operator and an error plot
for $\lambda=0,1,2$ resulting in $s=-2,0,2$ for the function
$f(\phi,\theta)=\cos(\theta)$. First one can see that for a
non-uniform probability measure the results for different values
of $\lambda$ differ quite a lot. Note that the function $f$ is
adapted to the cluster structure in the sense that it does not
change much in each cluster but changes very much in region of low
density. In the case of $s=2$ we can see that $\Delta_s f$ would
lead to a diffusion which would lead roughly to a kind of step
function which changes at the equator. The same is true for $s=0$
but the effect is much smaller than for $s=2$. In the case of
$s=-2$ we have a completely different behavior. $\Delta_s f$ has
now flipped its sign near to the equator so that the induced
diffusion process would try to smooth the function in the low
density region.

\section{Proof of the Main Result}\label{sec:continuum-limit}
In this section we will present the main results which were
sketched in Section \ref{sec:limit-sketch} together with the
proofs.  In Section \ref{sec:submanifolds-operators} we first
introduce some non-standard tools from differential geometry which
we will use later on. in particular,  it turns out that the so
called manifolds with boundary of bounded geometry are the natural
framework where one can still deal with non-compact manifolds in a
setting comparable to the compact case. After a proper statement
of the assumptions under which we prove the convergence results of
the graph Laplacian and a preliminary result about convolutions on
submanifolds which is of interest on its own, we then start with
the final proofs. The proof is basically divided into two parts,
the bias and the variance, where these terms are only
approximately valid. The reader not familiar with differential
geometry is encouraged to first read the appendix on basics of
differential geometry in order to be equipped with the necessary
background.

\subsection{Non-compact Submanifolds in $\R^d$ with Boundary}\label{sec:submanifolds-operators}
We prove the pointwise convergence for non-compact submanifolds.
Therefore we have to restrict the class of submanifolds since
manifolds with unbounded curvature do not allow reasonable
function spaces.

\begin{remark}
In the rest of this paper we use the Einstein summation convention
that is over indices occurring twice has to be summed. Note that
the definition of the curvature tensor differs between textbooks.
We use here the conventions regarding the definitions of curvature
etc. of \cite{Lee1997}.
\end{remark}

\subsubsection{Manifolds with Boundary of Bounded Geometry}\label{sec:bounded-geometry}
We will consider in general non-compact submanifolds with
boundary. In textbooks on Riemannian geometry one usually only
finds material for the case where the manifold has no boundary.
Also the analysis e.g. definition of Sobolev spaces on non-compact
Riemannian manifolds seems to be non-standard. We profit here very
much from the thesis and an accompanying paper of
\cite{Schick1996,Schick2001} which introduces manifolds with
boundary of bounded geometry. All material of this section is
taken from these articles. Naturally this plus of generality leads
also to a slightly larger technical overload. Nevertheless we
think that it is worth this effort since the class of manifolds
with boundary of bounded geometry includes almost any kind of
submanifold one could have in mind. Moreover, to our knowledge, it
is the most general setting where one can still introduce a notion
of Sobolev spaces with the usual properties.\\
Note that the boundary $\partial M$ is an isometric submanifold of
$M$ of dimension $m-1$. Therefore it has a second fundamental form
$\overline{\Pi}$ which should not be mixed up with the second
fundamental form $\Pi$ of $M$ which is with respect to the ambient
space $\R^d$. We denote by $\overline{\nabla}$ the connection and
by $\overline{R}$ the curvature of $\partial M$.
        \newif\iftobechos
        \tobechosfalse
        \iftobechos
 Moreover,  let $\nu$ be
the normal inward vector field at $\partial M$ and let $K$ be the
normal geodesic flow defined as $K: \partial M \times \R_+
\rightarrow M: (x',t) \rightarrow \exp_{x'}^M(t\nu_{x'})$. Then
the collar set $N(s)$ is defined as $N(s)\bydef K(\partial M
\times [0,s])$ for $s\geq 0$.
\begin{definition}[Manifold with boundary of bounded geometry]\label{def:manifold-bounded-geometry}
  Suppose $M$ is a manifold with boundary $\partial M$ (possibly empty). It is of bounded geometry if the
following holds:
\begin{itemize}
\item (N) Normal Collar: there exists $r_C>0$ so that the geodesic
collar
\[ \partial M \times [0,r_C) \rightarrow M: (x,t) \rightarrow \exp_x(t\nu_x)\]
is a diffeomorphism onto its image ($\nu_x$ is the inward normal
vector).
        \else
Moreover,  let $\nu$ be the normal inward vector field at
$\partial M$.

\begin{definition}[Manifold with boundary of bounded geometry]\label{def:manifold-bounded-geometry}
Let $M$ be a manifold with boundary $\partial M$ (possibly empty).
It is of bounded geometry if the following holds:
\begin{itemize}
\item (N) Normal Collar: there exists $r_C>0$ so that the normal
geodesic flow
\[ K: (x,t) \rightarrow \exp_x(t\nu_x)\]
is defined on $\partial M \times [0,r_C)$ and is a diffeomorphism
onto its image ($\nu_x$ is the inward normal vector). Let
$N(s)\bydef K(\partial M \times [0,s])$ be the collar set for $0
\le s\le r_C$.
        \fi
\item  (IC) The injectivity radius $\inj_{\partial M}$ of
$\partial M$ is positive. \item  (I) Injectivity radius of $M$:
There is $r_i>0$ so that if $r\leq r_i$ then for $x \in
M\backslash N(r)$
       the exponential map is a diffeomorphism on $B_M(0,r) \subset T_x M$ so that normal coordinates are
       defined on every ball $B_M(x,r)$ for $x \in M \backslash N(r)$.
\item  (B) Curvature bounds: For every $k \in \N$ there is $C_k$
so that $|\nabla^i R|\leq C_k$ and
        $\overline{\nabla}^i \overline{\Pi}\leq C_k$ for $0\leq i \leq k$, where $\nabla^i$ denotes the covariant derivative of order $i$.
\end{itemize}
\end{definition}
Note that $(B)$ imposes bounds on all orders of the derivatives of
the curvatures. One could also restrict the definition to the
order of derivatives needed for the goals one pursues. But this
would require even more notational effort, therefore we skip this.
in particular,  in \cite{Schick1996} it is argued that boundedness
of all derivatives
of the curvature is very close to the boundedness of the curvature alone.\\
The lower bound on the injectivity radius of $M$ and the bound on
the curvature are standard to define manifolds of bounded geometry
without boundary. Now the problem of the injectivity radius of $M$
is that at the boundary it somehow makes only partially sense
since $\inj_M(x)\rightarrow 0$ as $d(x,\partial M)\rightarrow 0$.
Therefore one replaces next to the boundary standard normal
coordinates with normal collar coordinates.
\begin{definition}[normal collar coordinates]
Let $M$ be a Riemannian manifold with\\ boundary $\partial M$. Fix
$x' \in \partial M$ and an orthonormal basis of $T_{x'}\partial M$
to identify $T_{x'}\partial M$ with $\R^{m-1}$. For $r_1,r_2>0$
sufficiently small (such that the following map is injective)
define normal collar coordinates,
\[ n_{x'}: B_{\R^{m-1}}(0,r_1)\times [0,r_2] \rightarrow M: (v,t) \rightarrow
\exp^M_{\exp^{\partial M}_{x'}(v)}(t\nu). \] The pair $(r_1,r_2)$
is called the width of the normal collar chart $n_{x'}$.
\end{definition}
The next proposition shows why manifolds of bounded geometry are
especially interesting.
\begin{proposition}[\cite{Schick2001}]\label{pro:metric-bounds}
Assume that conditions $(N),(IC),(I)$ of Definition
\ref{def:manifold-bounded-geometry} hold.
\begin{itemize}
\item(B1) There exist $0<R_1\leq r_\inj(\partial M)$, $0<R_2\leq
r_C$ and $0<R_3\leq r_i$ and constants $C_K>0$ for each $K \in \N$
such that whenever we have normal boundary coordinates of width
$(r_1,r_2)$ with $r_1\leq R_1$ and $r_2\leq R_2$ or normal
coordinates of radius $r_3\leq r_i$ then in these coordinates
\[ |D^\alpha g_{ij}| \leq C_K \quad \mathrm{and} \quad |D^\alpha g^{ij}| \leq C_K \quad \mathrm{for all}
\quad |\alpha|\leq K.\]
\end{itemize}
The condition $(B)$ in Definition
\ref{def:manifold-bounded-geometry} holds if and only if $(B1)$
holds. The constants $C_K$ can be chosen to depend only on
$r_i,r_C,\inj_{\partial M}$ and $C_k$.
\end{proposition}
Note that due to $g^{ij}g_{jk}=\delta^i_k$ one gets upper and
lower bounds on the operator norms of $g$ and $g^{-1}$,
respectively, which result in upper and lower bounds for
$\sqrt{\det g}$.
This implies that we have upper and lower bounds on the volume
form $dV(x)=\sqrt{\det g}\,dx$.
\begin{lemma}[\cite{Schick2001}]\label{le:volume-bounds}
Let $(M,g)$ be a Riemannian manifold with boundary of \\ bounded
geometry of dimension $m$. Then there exists $R_0>0$ and constants
$S_1>0$ and $S_2$ such that for all $x \in M$ and $r\leq R_0$ one
has
\[  S_1 r^m \leq \vol(B_M(x,r)) \leq S_2 r^m \]
\end{lemma}
Another important tool for analysis on manifolds are appropriate
function spaces. In order to define a Sobolev norm one first has
to fix a family of charts $U_i$ with $M \subset \cup_i U_i$ and
then define the Sobolev norm with respect to these charts. The
resulting norm will depend on the choice of the charts $U_i$.
Since in differential geometry the choice of the charts should not
matter, the natural question arises how the Sobolev norm
corresponding to a different choice of charts $V_i$ is related to
that for the choice $U_i$. In general, the Sobolev norms will not
be the same. However, if one assumes that the transition maps are
smooth and the manifold $M$ is compact then the resulting norms
will be equivalent and therefore define the same topology. Now if
one has a non-compact manifold this argumentation does not work
anymore. This problem is solved in general by defining the norm
with respect to a covering of $M$ by normal coordinate charts.
Then it can be shown that the change of coordinates between these
normal coordinate charts is well-behaved due to the bounded
geometry of $M$. In that way it is possible to establish a
well-defined notion of Sobolev spaces on manifolds with boundary
of bounded geometry in the sense that any norm defined with
respect to a different covering of $M$ by normal coordinate charts
is equivalent. Let $(U_i,\phi_i)_{i \in I}$ be a countable
covering of the submanifold $M$ with normal coordinate charts of
$M$, that is $M \subset \cup_{i \in I}U_i$, then:
\[ \norm{f}_{C^k(M)}= \max_{m\leq k} \sup_{i \in I} \sup_{x \in \phi_i(U_i)}
                   \big| D^m (f \circ \phi_i^{-1})(x)\big|.\]
In the following we will denote with $C^k(M)$ the space of
$C^k$-functions on $M$ together with the norm
$\norm{\cdot}_{C^k(M)}$.

\subsubsection{Intrinsic versus Extrinsic Properties}\label{sec:ex-in-comparison}
Most of the proofs for the continuous part will work with Taylor
expansions in normal coordinates. It is then of special interest
to have a connection between intrinsic and extrinsic distances.
Since the distance on $M$ is induced from $\R^d$, it is obvious
that one has $\norm{x-y}_{\R^d} \sim d_M(x,y)$ for all $x,y \in M$
which are sufficiently close. The next proposition proven by
 \citet*{SmoWeiWit2000}
provides an asymptotic expression of geometric quantities of the
submanifold $M$ in the neighborhood of a point $x \in M$.
Particularly, it gives a third-order approximation of the
intrinsic distance $d_M(x,y)$ in $M$ in terms of the extrinsic
distance in the ambient space $X$ which is in our case just the
Euclidean distance in $\R^d$.
\begin{proposition}\label{pro:normal-coordinates}
Let $i: M \rightarrow \R^d$ be an isometric embedding of the
smooth $m$-dimensional Riemannian manifold $M$ into $\R^d$. Let $x
\in M$ and $V$ be a neighborhood of $0$ in $\R^m$ and let $\Psi: V
\rightarrow U$ provide normal coordinates of a neighborhood $U$ of
$x$, that is $\Psi(0)=x$. Then for all $y \in V$:
\begin{equation*}
\norm{y}^2_{\R^m}=d^2_M(x,\Psi(y))=\norm{(i \circ \Psi)(y)-
i(x)}^2 +
\frac{1}{12}\norm{\Pi(\dot{\gamma},\dot{\gamma})}_{T_x\R^d}^2+
O(\norm{y}^5_{\R^m}), \nonumber
\end{equation*}
where $\Pi$ is the second fundamental form of $M$ and $\gamma$ the
unique geodesic from $x$ to $\Psi(y)$ such that $\dot{\gamma}=y^i
\partial_{y^i}$. The volume form $dV=\sqrt{\det \,g_{ij}(y)}\,dy$
of $M$ satisfies in normal coordinates,
\begin{equation*}
dV=\Big(1+\frac{1}{6}R_{iuvi}\,y^u y^v +
O(\norm{y}^3_{\R^m})\Big)dy, \nonumber
\end{equation*}
in particular
\[ (\Delta \sqrt{\det \,g_{ij}})(0)=-\frac{1}{3}R, \]
where $R$ is the scalar curvature (i.e.,
$R=g^{ik}g^{jl}R_{ijkl}$).
\end{proposition}
We would like to note that in \cite{SmoWeiWit2004} this
proposition was formulated for general ambient spaces $X$, that is
arbitrary Riemannian manifolds $X$. Using the more general form of
this proposition one could extend the results in this paper to
submanifolds of other ambient spaces $X$. However,  in order to
use the scheme one needs to know the geodesic distances in $X$,
which are usually not available for general Riemannian manifolds.
Nevertheless, for some special cases like the sphere,
one knows the geodesic distances. Submanifolds of the sphere could be of interest, for example in geophysics or astronomy.\\
The previous proposition is very helpful since it gives an
asymptotic expression of the geodesic distance $d_M(x,y)$ on $M$
in terms of the extrinsic Euclidean distance. The following lemma
is a non-asymptotic statement taken from \cite{BerEtAl2001} which
we present in a slightly different form. But first we establish a
connection between what they call the 'minimum radius of
curvature' and upper bounds on the extrinsic curvatures of $M$ and
$\partial M$. Let
\[  \Pi_{\max} = \sup_{x \in M}\, \sup_{v \in T_x M, \norm{v}=1}\norm{\Pi(v,v)}, \quad \quad \overline{\Pi}_{\max}= \sup_{x \in \partial M} \,\sup_{v \in T_x \partial M, \norm{v}=1}
   \norm{\overline{\Pi}(v,v)},\]
where $\overline{\Pi}$ is the second fundamental form of $\partial
M$ as a submanifold of $M$.
We set $\overline{\Pi}_{\max}=0$ if the boundary $\partial M$ is empty.\\
Using the relation between the acceleration in the ambient space
and the second fundamental form for unit-speed curves $\gamma$
with no acceleration in $M$ ($D_t \dot{\gamma}=0$) established in
section \ref{sec:second-fundamental-form}, we get for the
Euclidean acceleration of such a curve $\gamma$ in $\R^d$,
\[ \norm{\ddot{\gamma}}=\norm{\Pi(\dot{\gamma},\dot{\gamma})}. \]
Now if one has a non-empty boundary $\partial M$ it can happen
that a length-minimizing curve goes (partially) along the boundary
(imagine $\R^d$ with a ball at the origin cut out). Then the
segment $c$ along the boundary will be a geodesic of the
submanifold $\partial M$, see \cite{AleAle1981}, that is
$\overline{D}_t \dot{c}=\overline{\nabla}_{\dot{c}}\dot{c}=0$
where $\overline{\nabla}$ is the connection of $\partial M$
induced by $M$. However,  $c$ will not be a geodesic in $M$ (in
the sense of a curve with no acceleration) since by the
Gauss-Formula in Theorem \ref{th:gauss-formula},
\[ D_t \dot{c}= \overline{D}_t \dot{c} + \overline{\Pi}(\dot{c},\dot{c}) = \overline{\Pi}(\dot{c},\dot{c}).\]
Therefore, in general the upper bound on the Euclidean
acceleration of a length-minimizing curve $\gamma$ in $M$ is given
by,
\[ \norm{\ddot{\gamma}}=\norm{\overline{\Pi}(\dot{\gamma},\dot{\gamma})+\Pi(\dot{\gamma},\dot{\gamma})}
\leq \overline{\Pi}_{\max}+ \Pi_{\max}. \] Using this inequality,
one can derive a lower bound on the 'minimum radius of curvature'
$\rho$ defined in \cite{BerEtAl2001} as
$\rho=\inf\{1/\norm{\ddot{\gamma}}_{\R^d}\}$ where the infimum is
taken over all unit-speed geodesics $\gamma$ of $M$ (in the sense
of length-minimizing curves):
\[ \rho \geq \frac{1}{\overline{\Pi}_{\max}+ \Pi_{\max}}. \]
Finally we can formulate the Lemma from \cite{BerEtAl2001}.
\begin{lemma}
Let $x,y \in M$ with $d_M(x,y)\leq \pi \rho$. Then
\[ 2 \rho \sin(d_M(x,y)/(2\rho))\leq \norm{x-y}_{\R^d} \leq d_M(x,y). \]
\end{lemma}
Noting that $\sin(x)\geq x/2$ for $0\leq x \leq \pi/2$, we get as
an easier to handle corollary:
\begin{corollary}\label{co:comparison-ex-in}
Let $x,y \in M$ with $d_M(x,y)\leq \pi \rho$. Then
\[ \frac{1}{2} d_M(x,y) \leq \norm{x-y}_{\R^d} \leq d_M(x,y). \]
\end{corollary}
In the given form this corollary is quite useless since we only
have the Euclidean distances between points and therefore we have
no possibility to check the condition $d_M(x,y)\leq \pi \rho$. In
general small Euclidean distance does not imply small intrinsic
distance. Imagine a circle where one has cut out a very small
segment. Then the Euclidean distance between the two ends is very
small however the geodesic distance is very large. We show now
that under an additional assumption one can transform the above
corollary so that one can use it when one has only knowledge about
Euclidean distances.
\begin{lemma}\label{le:comparison-ex-in}
Let $M$ have a finite radius of curvature $\rho>0$. We further
assume that,
\[ \kappa \bydef \inf_{x \in M} \inf_{y \in M \backslash B_M(x,\pi \rho)}\norm{x-y}, \]
is non-zero. Then $B_{\R^d}(x,\kappa / 2)\cap M \subset
B_M(x,\kappa) \subset B_M(x,\pi\rho)$. Particularly, if $x,y \in
M$ and $\norm{x-y}\leq \kappa / 2$,
\[ \frac{1}{2} d_M(x,y) \leq \norm{x-y}_{\R^d} \leq d_M(x,y) \leq \kappa. \]
\end{lemma}
\begin{proof}
By definition $\kappa$ is at most the infimum of $\norm{x-y}$
where $y$ satisfies $d_M(x,y)=\pi\rho$. Therefore the set
$B_{\R^d}(x,\kappa / 2) \cap M$ is a subset of $B_M(x,\pi\rho)$.
The rest of the lemma then follows by Corollary
\ref{co:comparison-ex-in}. Figure \ref{fig:kappa} illustrates this
construction.
\end{proof}
\begin{figure}[h]
\vspace{-5mm}
$$\epsfig{file=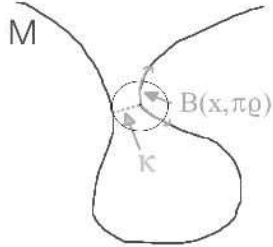,width=4cm}$$
\vspace{-10mm} \caption{$\kappa$ is the Euclidean distance of
$x\in M$ to $M\backslash B_M(x,\pi\rho)$. \label{fig:kappa}}
\end{figure}

\subsection{Notations and Assumptions}\label {sec:assumptions}
In general we work on complete non-compact manifolds with
boundary. Compared to a setting where one considers only compact
manifolds one needs a slightly larger technical overhead. However,
we will indicate how the technical assumptions simplify if one has
a compact submanifold with boundary or even
a compact manifold without boundary.\\
We impose the following assumptions on the manifold $M$:
\begin{assumptions}\label{manifold-assumptions}
\begin{enumerate}[(i)]
\item The map $i:M \rightarrow \R^d$ is a smooth embedding, \item
The manifold $M$ with the metric induced from $\R^d$ is a smooth
manifold with boundary of bounded geometry (possibly $\partial
M=\emptyset$), \item $M$ has bounded second fundamental form,
\item It holds $\kappa \bydef \inf_{x \in M} \inf_{y \in M
\backslash B_M(x,\pi \rho)}\norm{i(x)-i(y)}>0$, where $\rho$ is
the radius of curvature defined in Section
\ref{sec:ex-in-comparison},

\item For any $x \in M\backslash \partial M$,
$\delta(x):=\inf\limits_{y \in M \backslash
B_M(x,\frac{1}{3}\min\{\mathrm{inj}(x),\pi\rho\})}\norm{i(x)-i(y)}_{\R^d}>0,$
      where $\mathrm{inj}(x)$ is the injectivity radius\footnote{Note that the injectivity radius $\mathrm{inj}(x)$ is always positive.} at $x$
      and $\rho>0$ is the radius of curvature.
\end{enumerate}
\end{assumptions}
The first condition ensures that $i(M)$ is a smooth submanifold of
$\R^d$. Usually we do not distinguish between $i(M)$ and $M$. The
use of the abstract manifold $M$ as a starting point emphasizes
that there exists an $m$-dimensional smooth manifold $M$ or
roughly equivalent an $m$-dimensional smooth parameter space
underlying the data. The choice of the $d$ features determines
then the representation in $\R^d$. The choice of features
corresponds therefore to a specific choice of the inclusion map
$i$ since $i$ determines how $M$ is embedded into $\R^d$. This
means that another choice of features leads in general to a
different mapping $i$ but the initial abstract manifold $M$ is
always the same. However, in the second condition we assume that
the metric structure of $M$ is induced by $\R^d$ (which implies
that $i$ is trivially an isometric embedding). Therefore the
metric structure depends on the embedding $i$ or equivalently on
our choice of features.

The second condition ensures that $M$ is an isometric submanifold
of $\R^d$ which is well-behaved. As discussed in section
\ref{sec:bounded-geometry}, manifolds of bounded geometry are
 in general non-compact, complete Riemannian manifolds with boundary where one has uniform control over all intrinsic curvatures.
The uniform bounds on the curvature allow to do reasonable
analysis in this general setting. In particular, it allows us to
introduce the function spaces $C^k(M)$ with their associated norm.
It might be possible to prove pointwise results even without the
assumption of bounded geometry. But we think that the setting
studied here is already general enough to encompass all cases
encountered in practice.
The third condition ensures that $M$ also has well-behaved
extrinsic geometry and implies that the radius of curvature $\rho$
is lower bounded. Together with the fourth condition it enables us
to get global upper and lower bounds of the intrinsic distance on
$M$ in terms of the extrinsic distance in $\R^d$ and vice versa,
see Lemma \ref{le:comparison-ex-in}. The fourth condition is only
necessary in the case of non-compact submanifolds. It prevents the
manifold from self-approaching. More precisely it ensures that if
parts of $M$ are far away from $x$ in the geometry of $M$ they do
not come too close to $x$ in the geometry of $\R^d$. Assuming that
$i(M)$ is a submanifold, this assumption is already included
implicitly. However,  for non-compact submanifolds the
self-approaching could happen at infinity. Therefore we exclude it
explicitly. Moreover,  note that for submanifolds with boundary
one has $\mathrm{inj}(x)\rightarrow 0$ as $x$ approaches the
boundary\footnote{This is the reason why one replaces normal
coordinates in the neighborhood of the boundary with normal collar
coordinates.} $\partial M$. Therefore also $\delta(x)\rightarrow
0$ as $d(x,\partial M)\rightarrow 0$. However,  this behavior of
$\delta(x)$ at the boundary does not matter for the proof of
pointwise convergence in the interior of $M$. Note that if $M$ is
a smooth
and \emph{compact} manifold conditions (ii)-(v) hold automatically.\\
In order to emphasize the distinction between extrinsic and
intrinsic properties of the ma\-ni\-fold we always use the
slightly cumbersome notations $x \in M$ (intrinsic) and $i(x)\in
\R^d$ (extrinsic). The reader who is not familiar with Riemannian
geometry should keep in mind that locally, a submanifold of
dimension $m$ looks like $\R^m$. This becomes apparent if one uses
normal coordinates. Also the following dictionary between terms of
the manifold $M$ and the case when one has only an open set in
$\R^d$ ($i$ is then the identity mapping) might be useful.
\begin{center}
\begin{tabular}{|  c  |  c |}
\hline\multicolumn{1}{| c|}{{\centering{Manifold $M$}}}
&\multicolumn{1}{c|}{{\centering{open set in $\R^d$}}}\\
\hline
$g_{ij}$ , $\sqrt{\det g}$  & $\delta_{ij}$ , $1$ \\
natural volume element & Lebesgue measure\\
$\Delta_s$      &  $\Delta_s=\sum_{i=1}^d \frac{\partial^2}{\partial (z_i)^2} + \frac{s}{p}\sum_{i=1}^d \frac{\partial p}{\partial z^i}\frac{\partial}{\partial z^i}$\\
\hline
\end{tabular}
\end{center}
The kernel functions which are used to define the weights of the
graph are always functions of the squared norm in $\R^d$.
Furthermore, we make the following assumptions on the kernel
function $k$:
\begin{assumptions}\label{kernel-assumptions}
\begin{enumerate}[(i)]
\item $k:\R_+ \rightarrow \R$ is measurable, non-negative and
non-increasing on $\R^*_+$, \item $k \in C^2(\R^*_+)$, that is in
particular $k$, $\frac{\partial k}{\partial x}$ and
$\frac{\partial^2 k}{\partial x^2}$ are bounded, \item $k$,
$|\frac{\partial k}{\partial x}|$ and $|\frac{\partial^2
    k}{\partial x^2}|$ have exponential decay: $\exists c,\alpha,A \in
\R_+$ such that for any $t \geq A$, $f(t)\leq c e^{-\alpha t}$,
where $f(t)=\max\{k(t),|\frac{\partial k}{\partial
x}|(t),|\frac{\partial^2 k}{\partial x^2}|(t)\}$, \item $k(0)=0$.
\end{enumerate}
\end{assumptions}
The assumption that the kernel is non-increasing could be dropped,
however it makes the proof and the presentation easier. Moreover,
in practice the weights of the neighborhood graph which are
determined by $k$ are interpreted as similarities. Therefore the
usual choice is to take weights which decrease with increasing
distance. The fourth condition implies that the graph has no
loops\footnote{An edge from a vertex to itself is called a loop.}.
in particular,  the kernel is not continuous at the origin. All
results hold also without this condition. The advantage of this
condition is that some estimators become unbiased. Also let us
introduce the helpful notation,
$k_h(t)=\frac{1}{h^m}k\left(\frac{t}{h^2}\right)$ where we call
$h$ the bandwidth of the kernel. Moreover,  we define the
following two constants related to the kernel function $k$,
\begin{equation}\label{eq:kernel-constants}
C_1 = \int_{\R^m} k(\norm{y}^2) dy < \infty, \quad C_2 =
\int_{\R^m} k(\norm{y}^2) y_1^2 dy < \infty.
\end{equation}
We also have some assumptions on the probability measure $P$.
\begin{assumptions}\label{probability-assumptions}
\begin{enumerate}[(i)]
\item $P$ is absolutely continuous with respect to the natural
volume element $dV$ on $M$, \item the density $p$ fulfills: $p \in
C^3(M)$ and $p(x)>0, \forall \; x \in M\backslash \partial M$,
\item the sample $X_i,\; i=1,\ldots,n$ is drawn i.i.d. from $P$,
\end{enumerate}
\end{assumptions}
Note that condition (i) implies $P(\partial M)=0$, that is the
boundary $\partial M$ is a set of measure zero. We will call the
Assumptions \ref{manifold-assumptions} on the submanifold,
Assumptions \ref{kernel-assumptions} on the kernel function, and
Assumptions \ref{probability-assumptions} on the probability
measure $P$ together the \textbf{standard assumptions}.


\setlength{\myVSpace}{0.7cm} In the following table we summarize
the notation used in the proofs:
\begin{center}\label{tab:notation1}
\begin{tabular}{|  c  |  c |}
\hline
$k:\R_+\rightarrow\R_+$ & kernel function \\
$h>0$ & neighborhood/bandwidth parameter\\
$m\in\N$ & dimension of the submanifold $M$\\ 
$k_h(t)=\frac{1}{h^m}k\left(\frac{t}{h^2}\right)$ & scaled kernel function \\
$\lambda\in\R$ & reweighting parameter\\
\xstrut $d_{h,n}(x)=\frac{1}{n}\sum_{i=1}^n k_h(\norm{x-X_i}^2)$ &  degree function associated with $k$\\
\xstrut $\tk_{\lambda,h}(x,X_i)=\frac{k_h(\|x-X_i\|^2)}{[d_{h,n}(X_i)d_{h,n}(X_j)]^\lambda}$ & reweighted kernel\\
\xstrut $\td_{\lambda,h,n}(x)=\frac{1}{n}\sum_{i=1}^n \tk_{\lambda,h}(x,X_i)$ & degree function associated with $\tk_{\lambda,h}$\\
\xstrut $(\tA_{\lambda,h,n}f)(x)=\frac{1}{n}\sum_{i=1}^n \tilde{k}_{\lambda,h}(x,X_i)f(X_i)$ & empirical average operator $\tA_{\lambda,h,n}$\\
\hline
\xstrut $\Delrw_{\lambda,h,n}f =\frac{1}{h^2}\Big(f- \frac{1}{\td_{\lambda,h,n}} \tA_{\lambda,h,n}f\Big)$ & random walk graph Laplacian\\
\xstrut $\Delu_{\lambda,h,n}f = \frac{1}{h^2}\big(\td_{\lambda,h,n}f - \tA_{\lambda,h,n}f \big)$ & unnormalized graph Laplacian\\
\xstrut $\Deln_{\lambda,h,n}f= \frac{1}{h^2}\Big(f - \frac{1}{\sqrt{\td_{\lambda,h,n}}}\tA_{\lambda,h,n}\Big(\frac{f}{\sqrt{\td_{\lambda,h,n}}}\Big)\Big)$ & normalized graph Laplacian\\
\hline
\xstrut$C_1 = \int_{\R^m} k(\norm{y}_{\R^m}^2) dy$, $C_2 = \int_{\R^m} k(\norm{y}_{\R^m}^2) y_1^2 dy$ & characteristic constants of the kernel\\
\xstrut$p_h(x)=\Exp_Z k_h(\norm{x-Z}^2)$ & convolution of $p$ with $k_h$\\
\xstrut$(\tA_{\lambda,h} f)(x) =\Exp_Z \tilde{k}_{\lambda,h}(x,Z)f(Z)$ & average operator $\tA_{\lambda,h}$\\
\xstrut$\Delrw_{\lambda,h}$, $\Delu_{\lambda,h}$, $\Deln_{\lambda,h}$ & Laplacians associated with $\tA_{\lambda,h}$\\
\xstrut$\Delta_s=\frac{1}{p^s}\div(p^s\,\grad)
= \frac{1}{p^s}g^{ab}\nabla_a (p^s \nabla_b)$      &  $s$-th weighted Laplacian on $M$\\
\hline
\end{tabular}
\end{center}

\subsection{Asymptotics of Euclidean Convolutions on the Submanifold $M$}\label{sec:convolution}
The following proposition describes the asymptotic expression of
the convolution of a function $f$ on the submanifold $M$ with a
kernel function having the Euclidean distance $\norm{x-y}$ as its
argument with respect to the probability measure $P$ on $M$. This
result is interesting since it shows how the use of the Euclidean
distance introduces a curvature effect if one averages a function
locally. A similar result has been presented in \cite{CoiLaf2005}.
We define the density $p$ invariantly with respect to the natural
volume element and also explicitly give the second order curvature
terms. Our proof is similar to that of \citet*{SmoWeiWit2004}
where under stronger conditions a similar result was proven for
the Gaussian kernel. The more general setting and the use of
general kernel functions make the proof slightly more complicated.
In order to emphasize the distinction between extrinsic and
intrinsic properties of the manifold we will use the slightly
cumbersome notations $x \in M$ (intrinsic) and $i(x)\in \R^d$
(extrinsic).

\begin{proposition}\label{pro:averaging-op-manifold}
Let $M$ and $k$ satisfy Assumptions \ref{manifold-assumptions} and
\ref{kernel-assumptions}. Furthermore,  let $P$ have a density $p$
with respect to the natural volume element and $p \in C^3(M)$.
Then, for any $x \in M\backslash \partial M$, there exists an
$h_0(x)>0$ such that for all $h<h_0(x)$ and any $f\in C^3(M)$,
\begin{align}
 \int_{M} &k_h\big(\norm{i(x)-i(y)}^2_{\R^d}\big)f(y)p(y)\sqrt{\det g}\,dy \nonumber \\
=&C_1 p(x)f(x) + \frac{h^2}{2} \,C_2 \Big(p(x)f(x)S(x) +
(\Delta_M(pf))(x)\Big) + O(h^3), \nonumber
\end{align}
where $O(h^3)$ is a function depending on $x$, $\norm{f}_{C^3(M)}$
and $\norm{p}_{C^3(M)}$ and
\[ S(x)=\frac{1}{2}\Big[-R\big|_x + \frac{1}{2}\norm{\sum\nolimits_a \Pi(\partial_a,\partial_a)}_{T_{i(x)}\R^d}^2\Big], \]
where $R$ is the scalar curvature and $\Pi$ the second fundamental
form of $M$ .
\end{proposition}
The following Lemma is an application of Bernstein's inequality.
Together with the previous proposition it will be the main
ingredient for proving consistency statements for the graph
structure.
\begin{lemma} \label{le:variance}
Suppose the standard assumptions hold and let the kernel $k$ have
compact support on $[0,R_k^2]$. Define $b_1=\norm{k}_\infty
\norm{f}_\infty, \; b_2=K\norm{f}_\infty^2$ where $K$ is a
constant depending on $\norm{p}_\infty$, $\norm{k}_\infty$ and
$R_k$. Let $x \in M \backslash \partial M$ and $V_i\bydef
k_h(\norm{i(x)-i(X_i)}^2) f(X_i)$. Then for any bounded function
$f$,
\begin{align}
\Pr \Big(\Big|\frac{1}{n} \sum_{i=1}^n V_i - \Exp V \Big| >
\epsilon\Big) \leq 2 \exp\biggl(- \frac{n h^m\eps^2}{2 b_2 + 2b_1
\eps/3}\biggr). \nonumber
\end{align}
Let $W_i=k_h(\norm{i(x)-i(X_i)}^2)(f(x)-f(X_i))$. Then for $hR_k
\leq \kappa/2$ and $f \in C^1(M)$,
\begin{align}
\Pr \Big(\Big|\frac{1}{n} \sum_{i=1}^n W_i  - \Exp W \Big| >
h\,\epsilon \Big) \leq 2 \exp\bigg(- \frac{n h^m\eps^2}{2 b_2 +
2b_1 \eps/3}\bigg). \nonumber
\end{align}
\end{lemma}
\begin{proof}
Since by assumption $\kappa>0$, by Lemma
\ref{le:comparison-ex-in}, for any $x,y \in M$ with
$\norm{i(x)-i(y)}\leq \kappa/2$, we have $d_M(x,y)\leq
2\norm{i(x)-i(y)}$. This implies $\forall a\leq \kappa/2$,
$B_{\R^d}(x,a)\cap M \subset B_M(x,2a)$.\\
Let $W_i \bydef k_h(\norm{i(x)-i(X_i)}^2) f(X_i)$. We have
\[|W_i| \le \frac{\norm{k}_\infty}{h^m} \sup_{y \in B_{\R^d}(x,hR_k) \cap M} |f(y)|
        \leq \frac{\norm{k}_\infty}{h^m}\norm{f}_\infty \bydef \frac{b_1}{h^m}.\]
For the variance of $W$ we have two cases. First let $hR_k<s
\bydef \min\{\kappa/2,R_0/2\}$.  Then we get
\begin{align*}
\Var W  \leq \Exp_Z k^2_h(\norm{i(x)-i(Z)}^2)f^2(Z) \leq
\frac{\norm{k}_\infty}{h^m}\norm{f}_\infty^2 p_h(x)
        \leq D_2 \frac{\norm{k}_\infty}{h^m}\norm{f}_\infty^2
\end{align*}
where we have used Lemma \ref{le:averaged-density-bounds} in the
last step. Now consider $hR_k \geq s$, then
\begin{align*}
\Var W &\leq \frac{\norm{k}_\infty^2}{h^{2m}}\norm{f}_\infty^2
\leq \frac{R_k^m \norm{k}_\infty^2}{s^m \,h^{m}}\norm{f}_\infty^2
\end{align*}
Therefore we define $b_2= K \norm{f}^2_\infty$ with $K = R_k^m
\norm{k}^2_\infty \max\{2^m S_2 \norm{p}_\infty, s^{-m}\}$. By
Bernstein's inequality we finally get
    \begar
    \Pr\Big(\big|\frac{1}{n}\sum_{i=1}^{n} W_i - \dsE W \big| > \eps\Big)
    & \le & 2 e^{- \frac{n h^m\eps^2}{2 b_2 + 2b_1 \eps/3}}\\
    \endar
Both constants $b_2$ and $b_1$ are independent of $x$. For the
second part note that by Lemma \ref{le:comparison-ex-in} for $h
R_k \leq \kappa /2$, we have that $\norm{x-y}\leq h R_k$ implies
$d_M(x,y)\leq 2 \norm{x-y}\leq 2 hR_k$. in particular,  for all
$x,y \in M$ with $\norm{x-y} \leq h R_k$,
\[ |f(x)-f(y)|\leq \sup_{y \in M} \norm{\nabla f}_{T_y M} d_M(x,y) \leq 2 h R_k  \sup_{y \in M} \norm{\nabla f}_{T_y M}. \]
A similar reasoning as above leads then to the second statement.
\end{proof}
Note that $\Exp_Z k_h(\norm{i(x)-i(Z)}^2) f(Z)=\int\limits_{M}
k_h(\norm{i(x)-i(y)}^2) f(y) p(y) \sqrt{\det g}\,dy$.

\subsection{Pointwise Consistency of the Random Walk, Unnormalized and Normalized Graph Laplacian}\label{sec:laplacian-consistency}
The proof of the convergence result for the three graph Laplacians
is organized as follows. First we introduce the continuous
operators $\Delrw_{\lambda,h},\; \Deln_{\lambda,h}$ and
$\Delu_{\lambda,h}$. Then we derive the limit of the continuous
operators as $h\rightarrow 0$. This part of the proof is concerned
with the bias part since roughly $(\Delta_{\lambda,h}f)(x)$ can be
seen as the expectation of $\Delta_{\lambda,h,n}f(x)$. Second we
show that with high probability all extended graph Laplacians are
close to the corresponding continuous operators. This is the
variance part. Combining both results we arrive finally at the
desired consistency results.
\subsubsection{The Bias Part - Deviation of $\Delta_{\lambda,h}$ from its Limit}
The following continuous approximation of $\Delrw_{\lambda,h}$ was
similarly introduced in \cite{Lafon2004,CoiLaf2005}.
\begin{definition}[Kernel-based approximation of the Laplacian]
We introduce the \\following averaging operator $\tA_{\lambda,h}$
based on the reweighted kernel $\tilde{k}_{\lambda,h}$:
\begin{equation}
(\tA_{\lambda,h} f)(x) =\int_{M} \tilde{k}_{\lambda,h}(x,y)f(y)
p(y) \sqrt{\det g}\, dy,
\end{equation}
and with $\td_{\lambda,h}=(\tA_{\lambda,h}1)$ the following
continuous operators:
\begin{align*}
\mathrm{random \, walk:} && \Delrw_{\lambda,h}f   &\bydef
\frac{1}{h^2}\bigg(f-
\frac{1}{\td_{\lambda,h}}\tA_{\lambda,h}f\bigg)
                        =\frac{1}{h^2}\bigg(\frac{\tA_{\lambda,h}g}{\td_{\lambda,h}} \bigg)(x),
                      \\
\mathrm{unnormalized:} && \Delu_{\lambda,h}f  &\bydef \frac{1}{h^2}\Big(\td_{\lambda,h}f - \tA_{\lambda,h}f\Big)= \frac{1}{h^2} (\tA_{\lambda,h}g)(x),
\\
\mathrm{normalized:} && \Deln_{\lambda,h}f &\bydef \frac{1}{h^2\, \sqrt{\td_{\lambda,h}}}\Big(d_{\lambda,h}\frac{f}{\sqrt{\td_{\lambda,h}}}- \tA_{\lambda,h}\frac{f}{\sqrt{\td_{\lambda,h}}}\Big) 
                      = \tfrac{1}{h^2\,\sqrt{\td_{\lambda,h}(x)}}(\tA_{\lambda,h}g')(x),
\end{align*}
\end{definition}
where we have introduced again $g(y):=f(x)-f(y)$ and
$g'(y):=\frac{f(x)}{\sqrt{\td_{\lambda,h}(x)}}-\frac{f(y)}{\sqrt{\td_{\lambda,h}(y)}}$.
The definition of the normalized approximation
$\Delrw_{\lambda,h}$ can be justified by the alternative
definition of the Laplacian in $\R^d$ sometimes made in physics
textbooks:
\[ (\Delta f)(x)=\lim_{r \rightarrow 0} -\frac{1}{C_d \,r^2}\Big(f(x) - \frac{1}{\vol(B(x,r))}\int_{B(x,r)}f(y)dy \Big), \]
where $C_d$ is a constant depending on the dimension $d$.\\
Approximations of the Laplace-Beltrami operator based on averaging
with the Gaussian kernel in the case of a uniform probability
measure have been studied for compact submanifolds without
boundary by \citet*{SmoWeiWit2000,SmoWeiWit2004} and
\cite{Belkin2003}. Their result was then generalized by
\cite{Lafon2004} to general densities and to a wider class of
isotropic, positive definite kernels for compact submanifolds with
boundary. The proof given in \cite{Lafon2004} applies only to
compact hypersurfaces\footnote{A hypersurface is a submanifold of
codimension $1$.} in $\R^d$, a proof for the general case of
compact submanifolds with boundary using boundary conditions has
been presented in \cite{CoiLaf2005}. In this section we will prove
the pointwise convergence of the continuous approximation for
general submanifolds $M$ with boundary of bounded geometry with
the additional Assumptions \ref{manifold-assumptions}. This
includes the case where $M$ is not compact. Moreover, no
assumptions of positive definiteness of the kernel are made nor
any boundary condition on the function $f$ is imposed. Almost any
submanifold occurring in practice should be
covered in this very general setting.\\
For pointwise convergence in the interior of the manifold $M$
boundary conditions on $f$ are not necessary. However,  for
uniform convergence there is no way around them. Then the problem
lies not in the proof that the continuous approximation still
converges in the right way but in the transfer of the boundary
condition to the discrete graph. The main problem is that since we
have no information about $M$ apart from the random samples the
boundary will be hard to locate. Moreover,  since the boundary is
a set of measure zero, we will actually almost surely never sample
any point from the boundary. The rigorous treatment of the
approximation of the boundary respectively the boundary conditions
of a function on a randomly sampled graph
remains as an open problem.\\
Especially for dimensionality reduction the case of
low-dimensional submanifolds in $\R^d$ is important. Notably, the
analysis below also includes the case where due to noise the data
is only concentrated around a submanifold.
\begin{theorem}\label{th:convergence-bias-manifold}
Suppose the standard assumptions hold. Furthermore,  let $k$ be a
kernel with compact support on $[0,R_k^2]$. Let $\lambda \in \R$,
and $x \in M\backslash \partial M$. Then there exists an
$h_1(x)>0$ such that for all $h<h_1(x)$  and any $f\in C^3(M)$,
\begin{align}
(\Delrw_{\lambda,h} f)(x) =&-\frac{C_2}{2\, C_1}\Big((\Delta_M f)(x) + \frac{s}{p(x)}\inner{\nabla p,\nabla f}_{T_xM}\Big) + O(h)
=-\frac{C_2}{2\, C_1}(\Delta_s f)(x)+O(h), \nonumber
\end{align}
where $\Delta_M$ is the Laplace-Beltrami operator of $M$ and
$s=2(1-\lambda)$.
\end{theorem}
\begin{proof}
For sufficiently small $h$ we have $\overline{B_{\R^d}(x,2hR_k)
\cap M} \cap \partial M = \emptyset$. Moreover,  it can be
directly seen from the
  proof of Proposition \ref{pro:averaging-op-manifold} that the upper bound of the interval $[0,h_0(y)]$ for which the expansion holds depends
  continuously on $\delta(x)$ and $\epsilon(y)$, where $\epsilon(y)=\frac{1}{3}\min\{\pi\rho,\mathrm{inj}(y)\}$. Now $h_0(x)$ is continuous since
  $\mathrm{inj}(x)$ is continuous on compact subsets, see \cite{Klingenberg1982}[Prop.  2.1.10], and $\delta(x)$ is continuous since the injectivity radius is continuous.
  Therefore we conclude that since $h_0(y)$ is continuous on $\overline{B(x,2hR_k)\cap M}$ and $h_0(y)>0$, $h_1(x)=\inf_{y \in
  \overline{B_{\R^d}(x,2hR_k)\cap M}}h_0(y)>0$. Then for the interval $(0,h_1(x))$ the expansion of $p_h(y)$ holds uniformly over the whole set
  $B(x,2hR_k)\cap M$.
  That is, using the definition of
  $\tilde{k}$ as well as Proposition \ref{pro:averaging-op-manifold}
  and the expansion $\frac{1}{(a+h^2
    b)^\lambda}=\frac{1}{a^\lambda}-\lambda \frac{ h^2
    b}{a^{\lambda+1}} + O(h^4)$, we get for $h \in (0,h_1(x))$ that
\begin{align*}
 &\int_{M} \tilde{k}_{\lambda,h}\bigl(\norm{i(x)-i(y)}^2 \bigr)f(y) p(y) \sqrt{\det g} \,dy \nonumber \\
=&\hspace{-1mm}  \int\limits_{B_{\R^d}(x,hR_k)\cap M} \hspace{-6mm} \frac{k_h\bigl(\norm{i(x)-i(y)}^2\bigr)}{p^\lambda_h(x)}f(y) 
 \bigg[\frac{C_1 p(y)- \lambda/2 C_2 h^2(p(y)S+\Delta p)}{C_1^{\lambda+1} p(y)^\lambda} +O(h^3)\bigg]\,\sqrt{\det g} \,dy,
\end{align*}
where the $O(h^3)$-term is continuous on $B_{\R^d}(x,hR_k)$ and we
have introduced the abbreviation $S=\frac{1}{2}[-R +
\frac{1}{2}\norm{\sum_a
\Pi(\partial_a,\partial_a)}_{T_{i(x)}\R^d}^2]$. Using $f(y)=1$ we
get,
\begin{align*}
 \td_{\lambda,h}(x)=\hspace{-4mm}&\int\limits_{B_{\R^d}(x,hR_k)\cap M}\hspace{-7mm} \frac{k_h\bigl(\norm{i(x)-i(y)}^2\bigr)}{p^\lambda_h(x)} 
 \bigg[\frac{C_1 p(y) - \lambda/2 C_2 h^2(p(y)S+\Delta p)}{C_1^{\lambda+1} p(y)^\lambda} +O(h^3)\bigg] \,\sqrt{\det g}\, dy,
\end{align*}
as an estimate for $\td_{\lambda,h}(x)$. Now using Proposition
\ref{pro:averaging-op-manifold} again, we arrive at:
\begin{align*}
\Delrw_{\lambda,h}f
 =\frac{1}{h^2}\frac{\td_{\lambda,h} f - \tA_{\lambda,h}f}{\td_{\lambda,h}} 
= -\frac{C_2}{2\, C_1}\Big(\Delta_M f +
\frac{2(1-\lambda)}{p}\inner{\nabla p,\nabla f} \Big) + O(h),
\nonumber
\end{align*}
where all $O(h)$-terms are finite on $B_{\R^d}(x,hR_k)\cap M$
since $p$ is strictly positive.
\end{proof}
Note that the limit of $\Delrw_{\lambda,h}$ has the opposite sign
of $\Delta_{s}$. This is due to the fact that the Laplace-Beltrami
operator on manifolds is usually defined as a negative definite
operator (in analogy to the Laplace operator in $\R^d$), whereas
the graph Laplacian is positive definite. But this varies through
the literature, thus the reader should be aware of the sign
convention.
\begin{remark}
The assumption of compact support of the kernel $k$ is only
necessary in the case of non-compact manifolds $M$. For compact
manifolds a kernel with non-compact support, such as a Gaussian
kernel, would work, too. The reason for compact support of the
kernel comes from the fact that for non-compact manifolds there
exists no lower bound on a strictly positive density. This in turn
implies that one cannot upper bound the convolution with the
reweighted kernel if one does not impose additional assumptions on
the density. In practice the solution of graph-based methods for
large-scale problems is usually only possible for sparse
neighborhood graphs. Therefore the compactness assumption of the
kernel is quite realistic and does not exclude relevant cases.
\end{remark}
With the relations
\begin{align*}
 (\Delu_{\lambda,h,n}f)(x)  &=\td_{\lambda,h,n}(x)(\Delrw_{\lambda,h,n}f)(x)\\
 (\Deln_{\lambda,h,n}f)(x) &=
 1 / \sqrt{\td_{\lambda,h,n}(x)} \,\big(\Delu_{\lambda,h,n} \Big(f / \sqrt{\td_{\lambda,h,n}}\Big)\big)(x) 
\end{align*}
one can easily adapt the last lines of the previous proof to
derive the following corollary.
\begin{corollary}\label{co:convergence-bias-manifold-unnormalized}
Under the assumptions of Theorem
\ref{th:convergence-bias-manifold}. Let $\lambda \in \R$ and $x
\in M\backslash \partial M$. Then there exists an $h_1(x)>0$ such
that for all $h<h_1(x)$ and any $f\in C^3(M)$,
\begin{align*}
(\Delu_{\lambda,h} f)(x) &=-\frac{C_2}{2 C_1^{2\lambda}}\,p(x)^{1-2\lambda}(\Delta_{s} f)(x) + O(h), \quad \mathrm{where}\quad s=2(1-\lambda),\\
(\Deln_{\lambda,h} f)(x)&=-\frac{C_2}{2 C_1}
\,p(x)^{\frac{1}{2}-\lambda} \Delta_{s}
\bigg(\frac{f}{p^{\frac{1}{2}-\lambda}}\bigg)(x) + O(h).
\end{align*}
\end{corollary}
\subsubsection{The Variance Part - Deviation of $\Delta_{\lambda,h,n}$ from $\Delta_{\lambda,h}$}
Before we state the results for the general case with
data-dependent weights we now treat the case $\lambda=0$, that is
we have non-data-dependent weights. There the proof is
considerably simpler and much easier to follow. Moreover,  as
opposed to the general case here we get convergence in probability
under slightly weaker conditions than almost sure convergence.
Since this does not hold for the normalized graph Laplacian in
that case we will  only provide the general proof.
\begin{theorem}[Weak and strong pointwise consistency for $\lambda=0$]
\label{th:pointwise-consistency-normal} Suppose the standard
assumptions hold. Furthermore,  let $k$ be a kernel with compact
support on $[0,R_k^2]$. Let $x \in M \backslash \partial M$ and $f
\in C^3(M)$. Then if $h \rightarrow 0$ and $nh^{m+2} \rightarrow
\infty$,
\begin{align*}
\lim_{n \rightarrow \infty} \,(\Delrw_{0,h,n} f)(x) &=-\frac{C_2}{2\,C_1}\,(\Delta_{2} f)(x)  &&\textrm{in probability},\\
\lim_{n \rightarrow \infty} \,(\Delu_{0,h,n}
f)(x)&=-\frac{C_2}{2}\,p(x)(\Delta_{2} f)(x)   &&\textrm{in
probability}.
\end{align*}
If even $nh^{m+2}/\log n  \rightarrow \infty$, then almost sure
convergence holds.
\end{theorem}
\begin{proof}
We give the proof for $\Delrw_{0,h,n}$. The proof for
$\Delu_{0,h,n}$ can be directly derived with the second statement
of Lemma \ref{le:variance} for the variance term together with
Corollary \ref{co:convergence-bias-manifold-unnormalized} for the
bias term. Similar to the proof for the Nadaraya-Watson regression
estimate of \cite{GreKryPaw1984}, we rewrite the estimator
$\Delrw_{0,h,n}f$ in the following form
\begin{equation}
(\Delrw_{0,h,n}f)(x)=\frac{1}{h^2}\bigg[\frac{(C_{0,h} \,f)(x) +
B_{1n}}{1+ B_{2n}}\bigg],
\end{equation}
where
\begin{align}
(C_{0,h}\,f)(x)&=\frac{\Exp_Z k_h(\norm{i(x)-i(Z)}^2)g(Z)}{\Exp_Z k_h(\norm{i(x)-i(Z)}^2)} ,\nonumber \\
B_{1n}&=\frac{\frac{1}{n}\sum_{j=1}^n k_h(\norm{i(x)-i(X_j)}^2)g(X_j) - \Exp_Z k_h(\norm{i(x)-i(Z)}^2)g(Z)}{\Exp_Z k_h(\norm{i(x)-i(Z)}^2)}, \nonumber\\
B_{2n}&=\frac{\frac{1}{n}\sum_{j=1}^n k_h(\norm{i(x)-i(X_j)}^2) -
\Exp_Z k_h(\norm{i(x)-i(Z)}^2)}{\Exp_Z k_h(\norm{i(x)-i(Z)}^2)},
\nonumber
\end{align}
with $g(y):=f(x)-f(y)$. In Theorem
\ref{th:convergence-bias-manifold} we have shown that for $x \in M
\backslash \partial M$,
\begin{equation}
\lim_{h\rightarrow 0}(\Delrw_{0,h}f)(x)=\lim_{h \rightarrow
0}\frac{1}{h^2}(C_{0,h}\,g)(x)
    =-\frac{ C_2}{2\,C_1}(\Delta_{2} f)(x).
\end{equation}
Using the lower bound of $p_h(x)=\Exp_Z k_h(\norm{i(x)-i(Z)}^2)$
derived in Lemma \ref{le:averaged-density-bounds} we can for $hR_k
\leq \kappa /2$ directly apply Lemma \ref{le:variance}. Thus there
exist constants $d_1$ and $d_2$ such that
\begin{align}{\footnotesize}
 \Pr(\,\left|B_{1n}\right|\geq h^2 t\,)&\leq 2\,\exp\biggl(-\frac{nh^{m+2}\,t^2}{2\norm{k}_\infty(d_2 + t \,h\,d_1/3)}\biggr). \nonumber
\end{align}
The same analysis can be done for $B_{2n}$. This shows convergence
in probability. Complete convergence (which implies almost sure
convergence) can be shown by proving for all $t>0$ the convergence
of the series $\sum_{n=0}^\infty \Pr\left(|B_{1n}|\geq
h^2t\right)<\infty$. A sufficient condition for that is $nh^{m+2}
/ \log n \rightarrow \infty$ as $n \rightarrow \infty$.
\end{proof}
The weak pointwise consistency of the unnormalized graph Laplacian
for compact submanifolds with the uniform probability measure
using the Gaussian kernel for the weights and $\lambda=0$ was
proven by \cite{BelNiy2005}. A more general result appeared
independently in \citep{HeiAudLux2005}. We prove here the limits
of all three graph Laplacians for general submanifolds with
boundary of bounded
geometry, general probability measures $P$, and general kernel functions $k$ as stated in our standard assumptions.\\
The rest of this section is devoted to the general case $\lambda
\neq 0$. We show that with high probability the extended graph
Laplacians $\Delta_{\lambda,h,n}$ are pointwise close to the
continuous operators $\Delta_{\lambda,h}$ when applied to a
function $f \in C^3(M)$. The following proposition is helpful.
\begin{proposition}\label{pro:limit-average-general}
Suppose the standard assumptions hold. Furthermore,  let $k$ be a
kernel with compact support on $[0,R_k^2]$. Fix $\lambda \in \R$
and let $x \in M\backslash \partial M$, $f\in C^3(M)$ and define
$g(y):=f(x)-f(y)$. Then there exists a constant $C$ such that for
any $\frac{2\norm{k}_\infty}{nh^m}<\epsilon<1/C$, $0<h<
\frac{\kappa}{2 \, R_k}$, the following events hold with
probability at least $1-C \,n \,e^\frac{-nh^m\epsilon^2}{C}$,
\[
|(\tA_{\lambda,h,n}g)(x) - (\tA_{\lambda,h}g)(x)| \leq \epsilon\,
h, \quad  \quad
 |\td_{\lambda,h,n}(x) - \td_{\lambda,h}(x)| \leq \epsilon.
\]
\end{proposition}
\begin{proof}
The idea of this proof is to show that several empirical
quantities which can be expressed as a sum of i.i.d. random
variables are close to their expectation. Then one can deduce that
also $(\tA_{\lambda,h,n}g)(x)$ will be close to
$(\tA_{\lambda,h}g)(x)$. The proof for $\td_{\lambda,h,n}$ can
then be easily adapted from the following. We consider here only
$\lambda\geq 0$, the proof for $\lambda<0$ is even simpler.
Consider the event $\E$ for which one has
\begin{small-lar}{}
    \text{for any } j\in\{1,\dots,n\}, \;
        \big| d_{h,n}(X_j) - p_h(X_j) \big| \le \eps\\
    \big| d_{h,n}(x) - p_h(x) \big| \le \eps\\
    \Big| \frac{1}{n}\sum\limits_{j=1}^n \frac{k_h(\norm{i(x)-i(X_j)}^2)|g(X_j)|}{[p_h(x)\,p_h(X_j)]^{\lambda}}
        - \int\limits_{M} k_h(\norm{i(x)-i(y)}^2)|g(y)| \frac{p(y)}{[p_h(x)\,p_h(y)]^{\lambda}}\sqrt{\det g}\,dy \Big| \le h\, \eps
\end{small-lar}
We will now prove that for sufficiently large $C$ the event $\E$
holds with probability at least $1-Cn e^{-\frac{nh^m\eps^2}{C}}$.
For the second assertion defining $\E$, we use Lemma
\ref{le:variance}
\[ \Pr(\,|d_{h,n}(x)-p_h(x)| > \epsilon\,) \leq 2 \,\exp\left(- \frac{n h^m\eps^2}{2 b_2 + 2b_1 \eps/3}\right),\]
where $b_1$ and $b_2$ are constants depending on the kernel $k$
and $p$. For the first term in the event $\E$ remember that
$k(0)=0$. We get for $\frac{\norm{k}_\infty}{nh^m} < \eps/2$ and
$1\leq j\leq n$,
\begin{small-ar}{}
    \Pr\Big( \,\big| \frac{1}{n} \sum_{i=1}^n k_h(\norm{i(X_j)-i(X_i)}^2)
        - p_h(x) \big| > \eps \, \Big| X_j \Big )
        \leq 2\exp\Bigl(- \frac{(n-1)h^m\eps^2}{8 b_2 + 4 b_1 \eps/3}\Bigr).
\end{small-ar}
This follows by
\begin{align*}
\Big | \frac{1}{n}\sum_{i=1}^n k_h(\norm{i(X_j)-i(X_i)}^2) -
p_h(X_j)\Big | \leq \Big |\frac{1}{n(n-1)}\sum_{i=1}^n
k_h(\norm{i(X_j)-i(X_i)}^2) \Big |\\ +
     \Big |\frac{1}{n-1}\sum_{i\neq j} k_h(\norm{i(X_j)-i(X_i)}^2) - p_h(X_j) \Big |
\end{align*}
where the first term is upper bounded by
$\frac{\norm{k}_\infty}{nh^m}$. First integrating wrt to the law
of $X_j$ (the right hand side of the bound is independent of
$X_j$) and then using a union bound, we get
\begin{small-ar}{}
    \Pr\Big( \text{for any } j\in\{1,\dots,n\}, \;
        \big| d_{h,n}(X_j)-  p_h(X_j) \big| \le \eps \Big)
        > 1 - 2 n \exp\left(- \frac{(n-1)h^m\eps^2}{8 b_2 + 4b_1 \eps/3}\right).
\end{small-ar}
Noting that $\frac{1}{p_h(x) p_h(y)}$ is upper bounded by Lemma
\ref{le:averaged-density-bounds} we get by Lemma \ref{le:variance}
for $h R_k\leq \kappa /2$ a Bernstein type bound for the
probability of the third event in $\E$. Finally, combining all
these results, we obtain that there exists a constant $C$ such
that for $h \leq \frac{\kappa}{2 \, R_k}$ and $\frac{2
\|k\|_\infty}{nh^m} \le \eps \le 1$, the event\footnote{The upper
bound on $\eps$ is here not necessary but allows to write the
bound more compactly.}  $\E$ holds with probability at least $1 -
C n e^{-\frac{nh^m \eps^2}{C}}$. Let us define \lbegar
    \B & \bydef & \int_{M} k_h(\norm{i(x)-i(y)}^2)(f(x)-f(y)) [p_h(x)\,p_h(y)]^{-\lam} p(y) \sqrt{\det g} dy\\
    \hat{\B} & \bydef & \frac{1}{n} \sum_{j=1}^n k_h(\norm{i(x)-i(X_j)}^2)(f(x)-f(X_j))\big[d_{h,n}(x)\,d_{h,n}(X_j)\big]^{-\lambda}\\
\rendar then $(\tA_{\lambda,h,n}g)(x)=\hat{\B}$ and
$(\tA_{\lambda,h}g)(x)=\B$. Let us now work only on the event
$\E$. By Lemma \ref{le:averaged-density-bounds} for any $y\in
B_{\R^d}(x,hR_k)\cap M$ there exist constants $D_1,D_2$ such that
$0<D_1 \le p_h(y) \le D_2$. Using the first order Taylor formula
of $[x\mapsto x^{-\lambda}]$, we obtain that for any $\lambda \ge
0$ and $a,b > \beta$, $\big| a^{-\lambda} - b^{-\lambda} \big| \le
\lambda \beta^{-\lambda-1} |a-b|$. So we can write for $\eps <
D_1/2$, \begar
 \Big| \frac{1}{\big(d_{h,n}(x)\, d_{h,n}(X_j)\big)^{\lambda}} - \frac{1}{\big(p_h(x)\,p_h(X_j)\big)^{\lambda}} \Big|
 &\leq \lambda (D_1 - \eps)^{-2\lambda-2}|d_{h,n}(x)d_{h,n}(X_j) - p_h(x)p_h(X_j)| \\
 &\leq 2\,\lambda  (D_1 - \eps)^{-2\lambda-2} (D_2 +\epsilon)\epsilon:=C\, \epsilon.
\endar
Noting that for $h R_k\leq \kappa/2$ by Lemma
\ref{le:comparison-ex-in}, $d_M(x,y)\leq 2 h R_k$, $\forall y \in
B_{\R^d}(x,hR_k)\cap M$,
    \begar
    \big| \hat{\B} - \B \big|
        & \le & \Big| \frac{1}{n}\sum_{j=1}^n k_h(\norm{i(x)-i(X_j)}^2)|f(x)-f(X_j)| \,C \,\epsilon \\
    & &
        + \Big| \frac{1}{n}\sum_{j=1}^n k_h(\norm{i(x)-i(X_j)}^2)(f(x)-f(X_j))[p_h(x)\,p_h(X_j)]^{-\lambda}
        - \B \Big|\\
    & \leq & 2 C \norm{k}_\infty R_k \sup_{y \in M} \norm{\nabla f}_{T_y M}\, h \, \epsilon + h \, \epsilon
   \endar
We have proven that there exists a constant $C>1$ such that for
any $0<h< \frac{\kappa}{2 \, R_k}$ and $\frac{2
\|k\|_\infty}{nh^m} < \eps < 1/C$,
\[ \left|(\tA_{\lambda,h,n}g)(x)-(\tA_{\lambda,h}g)(x)\right|\leq C''' h \,\eps , \]
with probability at least $1-Cn e^{-\frac{nh^m\eps^2}{C}}$.
\end{proof}
This leads us to our first main result for the random walk and the
unnormalized graph Laplacian.
\begin{theorem}[Pointwise consistency of $\Delrw_{\lambda,h,n}$ and $\Delu_{\lambda,h,n}$]
\label{th:limit-laplacian} Suppose the standard assumptions hold.
Furthermore,  let $k$ be a kernel with compact support on
$[0,R_k^2]$. Let $x \in M\backslash \partial M$, $\lambda \in \R$.
Then  for any $f \in C^3(M)$ there exists a constant $C$ such that
for any $\frac{2\norm{k}_\infty}{nh^{m+1}}<\epsilon<1/C$,
$0<h<h_{\max}$ with probability at least $1-C \,n
\,e^\frac{-nh^{m+2}\epsilon^2}{C}$,
\begin{align*}
|(\Delrw_{\lambda,h,n}f)(x) - (\Delrw_{\lambda,h}f)(x)|  \, &\leq \,\epsilon,\\
|(\Delu_{\lambda,h,n}f)(x) - (\Delu_{\lambda,h}f)(x)| \,&\leq
\,\epsilon.
\end{align*}
Define $s=2(1-\lambda)$.  Then if $h \rightarrow 0$ and $nh^{m+2}
/ \log n \rightarrow \infty$,
\begin{align*}
\lim_{n \rightarrow \infty} \,(\Delrw_{\lambda,h,n}f)(x)
&= -\frac{C_2}{2\,C_1}(\Delta_{s} f)(x) &\textrm{almost surely},\\
\lim_{n \rightarrow \infty} \,(\Delu_{\lambda,h,n}f)(x) &=
-\frac{C_2}{2\,C_1^{2\lambda}}\,p(x)^{1-2\lambda}\,(\Delta_{s}
f)(x) &\textrm{almost surely}.
\end{align*}
in particular,  under the above conditions,
\begin{align*}
\Big| (\Delrw_{\lambda,h,n}f)(x) \, - \,
\big[-\frac{C_2}{2\,C_1}(\Delta_{s} f)(x)\big]\Big|
&=O(h) + O\left(\sqrt{\tfrac{\log n}{nh^{m+2}}}\right) & \textrm{a.s.}\,,\\
\Big| (\Delu_{\lambda,h,n}f)(x)\, -\,
\big[-\frac{C_2}{2\,C_1^{2\lambda}}\,p(x)^{1-2\lambda}\,(\Delta_{s}
f)(x)\big]\Big| &=O(h) + O\left(\sqrt{\tfrac{\log
n}{nh^{m+2}}}\right) & \textrm{a.s.}\,.
\end{align*}
The optimal rate for $h(n)$ is $h=O((\log n/ n)^\frac{1}{m+4})$.
\end{theorem}
\begin{proof}
In Equation \ref{eq:graph-laplacian-final} it was shown that
\[(\Delrw_{\lambda,h,n}f)(x) = \frac{1}{h^2}\left(\tfrac{\tA_{\lambda,h,n}g}{\td_{\lambda,h,n}} \right)(x), \quad
 (\Delu_{\lambda,h,n}f)(x)   = \frac{1}{h^2} (\tA_{\lambda,h,n}g)(x),\]
where $g(y):=f(x)-f(y)$. Since $f$ is Lipschitz we can directly
apply Proposition \ref{pro:limit-average-general} so that for the
unnormalized we get with probability $1-C \,n
\,e^\frac{-nh^{m+2}\epsilon^2}{C}$,
\[ |(\Delu_{\lambda,h,n}f)(x) - (\Delu_{\lambda,h}f)(x)| \leq \epsilon.\]
For the random walk Laplacian $\Delrw_{\lambda,h,n}$ we work on
the event where $|\td_{\lambda,h,n}-\td_{\lambda,h}|\leq h
\epsilon$, where $\epsilon \leq \frac{1}{2h}\td_{\lambda,h}$. This
holds by Proposition \ref{pro:limit-average-general} with
probability $1-C \,n \,e^\frac{-nh^{m+2}\epsilon^2}{C}$. Moreover,
note that by Lemmas \ref{le:comparison-ex-in} and
\ref{le:volume-bounds} for $hR_k \leq \min\{\kappa/2,R_0\}$, we
have
 \[ \big|\tA_{\lambda,h}g\big|\leq \frac{2^m R_k^m S_2}{D_1^{2\lambda}}\norm{p}_\infty\norm{k}_\infty\, 2\, L(f)\,h \,R_k.\]
Using Proposition \ref{pro:limit-average-general} for
$\tA_{\lambda,h,n}g$ and the bounds of $p_h(x)$ from Lemma
\ref{le:averaged-density-bounds},
\begin{align*}
\big|&(\Delrw_{\lambda,h,n}f)(x)-(\Delrw_{\lambda,h}f)(x)\big|
=\frac{1}{h^2}\left|\tfrac{(\tA_{\lambda,h,n}g)(x)}{\td_{\lambda,h,n}(x)}
                      -\tfrac{(\tA_{\lambda,h}g)(x)}{\td_{\lambda,h}(x)}\right|\\
&\leq \frac{1}{h^2}\biggl(\dfrac{|(\tA_{\lambda,h,n}g)(x)\, -
\,(\tA_{\lambda,h}g)(x)|}{\td_{\lambda,h,n}(x)} +
     (\tA_{\lambda,h}g)(x)\, \dfrac{|\td_{\lambda,h,n}(x)\, - \, \td_{\lambda,h}(x)|}{\td_{\lambda,h,n}(x)\td_{\lambda,h}(x)}\biggr)\\
&\leq \frac{2\,D_2^{2\lambda}}{D_1}\,\epsilon + \frac{2^m R_k^m
S_2}{D_1^{2\lambda}}\norm{p}_\infty \norm{k}_\infty\, 2
\,L(f)\,R_k \,\epsilon:=C\,\epsilon,
\end{align*}
with probability $1-C \,n \,e^\frac{-nh^{m+2}\epsilon^2}{C}$. By
Theorem \ref{th:convergence-bias-manifold} and
\ref{co:convergence-bias-manifold-unnormalized} we have for
$s=2(1-\lambda)$,
\begin{align*}
\Big|(\Delrw_{\lambda,h}f)(x)\,-\,   \Big[-\frac{C_2}{2\,C_1}(\Delta_s f)(x)\Big]\Big|\leq C\, h,\\
\Big|(\Delu_{\lambda,h}f)(x)\, -\,
\Big[-\frac{C_2}{2\,C_1^{2\lambda}}\,p(x)^{1-2\lambda}\,(\Delta_{s}
f)(x)\Big]\Big|\leq C\,h.
\end{align*}
Combining both results together with the Borel-Cantelli-Lemma
yields almost sure convergence. The optimal rate for $h(n)$
follows by equating both order terms.
\end{proof}
Using the relationship between the unnormalized and the normalized
Laplacian the pointwise consistency can be easily derived.
However,  the conditions for convergence are slightly stronger
since the Laplacian is applied to the function
$f/\sqrt{\td_{\lambda,h,n}}$.
\begin{theorem}[Pointwise consistency of $\Deln_{\lambda,h,n}$]\label{th:limit-laplacian-normalized}
Suppose that the standard assumptions hold. Furthermore,  let $k$
be a kernel with compact support on $[0,R_k^2]$. Let $x \in
M\backslash \partial M$, $\lambda \in \R$. Then for any $f \in
C^3(M)$  there exists a constant $C$ such that for any
$\frac{2\norm{k}_\infty}{nh^{m+2}}<\epsilon<1/C$, $0<h<h_{\max}$
with probability at least $1-C
\,n^2\,e^\frac{-nh^{m+4}\epsilon^2}{C}$,
\[    \Big|(\Deln_{\lambda,h,n}f)(x) - (\Deln_{\lambda,h}f)(x)\Big| \leq \epsilon. \]
Define $s=2(1-\lambda)$.  Then if $h \rightarrow 0$ and $nh^{m+4}
/ \log n \rightarrow \infty$,
\[ \lim_{n \rightarrow \infty}\, (\Deln_{\lambda,h,n}f)(x)=-p(x)^{\frac{1}{2}-\lambda} \frac{C_2}{2 C_1} \Delta_{s} \biggl(\frac{f}{p^{\frac{1}{2}-\lambda}}\biggr)(x) \quad \textrm{almost surely}. \]
\end{theorem}
\begin{proof}
We reduce the case of $\Deln_{\lambda,h,n}$ to the case of
$\Delu_{\lambda,h,n}$. We work on the event where
\[  |\td_{\lambda,h,n}(x) - \td_{\lambda,h}(x)| \leq h^2\,\epsilon, \quad
\quad |\td_{\lambda,h,n}(X_i) - \td_{\lambda,h}(X_i)| \leq h^2\,
\epsilon, \quad \forall\, i=1,\ldots,n\] From Proposition
\ref{pro:limit-average-general} we know that this holds with
probability at least $1-C \,n^2
\,e^\frac{-nh^{m+4}\epsilon^2}{C}$. Working on this event we get
by a similar argumentation as in the proof of Theorem
\ref{th:limit-laplacian} that there exists a constant $C'$ such
that
\begin{align*}
 \Big|(\Deln_{\lambda,h,n}f)(x)- \tfrac{1}{\td_{\lambda,h}(x)}\Big(\Delu_{\lambda,h,n}\tfrac{f}{\sqrt{\td_{\lambda,h}}}\Big)(x)\Big|
=\frac{1}{h^2} \Big| \td_{\lambda,h,n}(x)f(x) \left[ \tfrac{1}{\td_{\lambda,h}(x)}-\tfrac{1}{\td_{\lambda,h,n}(x)}\right]\\
   + \sum_{i=1}^n \tilde{k}_{\lambda,h}(x,X_i)f(X_i)\left[\tfrac{1}{\sqrt{\td_{\lambda,h}(x)\td_{\lambda,h}(X_i)}}-
       \tfrac{1}{\sqrt{\td_{\lambda,h,n}(x)\td_{\lambda,h,n}(X_i)}}\right] \Big| \leq C' \epsilon.
\end{align*}
Noting that $\frac{f}{\td_{\lambda,h}}$ is Lipschitz since $f$ and
$\td_{\lambda,h}$ are Lipschitz and upper and lower bounded, on $M
\cap B_{\R^d}(x,h\,R_k)$ one can apply Theorem
\ref{th:limit-laplacian} to derive the first statement. The second
statement follows by Corollary
\ref{co:convergence-bias-manifold-unnormalized}.
\end{proof}

\vspace{0.3cm}
\noindent{\bf Acknowledgements}\\
We would like to thank Olaf Wittich and Bernhard Sch{\"o}lkopf for
helpful comments.

\begin{appendix}
\section{Basic Concepts of Differential Geometry} \label{sec:diffgeometry}
In this section we introduce the necessary basics of differential
geometry, in particular normal coordinates and submanifolds in
$\R^d$, used in this paper. Note that the definition of the
Riemann curvature tensor varies across textbooks which can result
in sign-errors. Throughout the paper we use the convention of
\cite{Lee1997}.
\subsection{Basics}
\begin{definition}
A $d$-dimensional \textbf{manifold $X$ with boundary} is a
topological (Hausdorff) space such that every point has a
neighborhood homeomorphic to an open subset of
$\HS^d=\{(x^1,\ldots,x^d) \in \R^d \big| x_1 \geq 0\}$. A
\textbf{chart} (or local coordinate system) $(U,\phi)$ of a
manifold $X$ is an open set $U \subset X$ together with a
homeomorphism $\phi: U \rightarrow V$ of $U$ onto an open subset
$V \subset \HS^d$. The coordinates $(x^1,\ldots,x^d)$ of $\phi(x)$
are called the coordinates of $x$ in the chart $(U,\phi)$. A
\textbf{$C^r$-atlas} $\A$ is a collection of charts $$\A
\triangleq \cup \{(U_\alpha,\phi_\alpha),\alpha\in I\},$$ where
$I$ is an index set, such that $X=\cup_{\alpha\in I} U_\alpha$ and
for any $\alpha,\beta\in I$ the corresponding \textbf{transition
map}
$$\phi_\beta \circ \phi_\alpha^{-1} \big|_{\phi_\alpha(U_\alpha \cap U_\beta)}: \phi(U_\alpha \cap U_\beta) \rightarrow \HS^d $$
is r-times continuously differentiable. A \textbf{smooth manifold
with boundary} is a manifold with boundary with a
$C^\infty$-atlas.
\end{definition}
For more technical details behind the definition of a manifold
with boundary we refer to \cite{Lee2003}. Note that the boundary
$\partial M$ of $M$ is a $(d-1)$-dimensional manifold without
boundary. In textbooks one often only finds the definition of a
manifold without boundary which can be easily recovered from the
above definition by replacing $\HS^d$ with $\R^d$. The interior
$M\backslash\partial M$ of the manifold $M$ is a manifold without
boundary.
\begin{definition}
A subset $M$ of a $d$-dimensional manifold $X$ is a
$m$-dimensional \textbf{submanifold $M$ with boundary} if every
point $x \in M$ is in the domain of a chart $(U,\phi)$ of $X$ such
that
\[ \phi: U \cap M \rightarrow \HS^m \times {a}, \quad \phi(x)=(x^1,\ldots,x^m,a^1,\ldots,a^{d-m})\]
where $a$ is a fixed element in $\R^{d-m}$. $X$ is called the
\textbf{ambient space} of $M$.
\end{definition}
This definition excludes irregular cases like intersecting
submanifolds or self-approaching submanifolds. In the following it
is more appropriate to take the following point of view. Let $M$
be an $m$-dimensional manifold. The smooth mapping $i:M
\rightarrow X$ is said to be an immersion if $i$ is differentiable
and the differential of $i$ has rank $m$ everywhere. An injective
immersion is called embedding if it is an homeomorphism onto its
image. In this case $i(M)$ is a submanifold of $X$. If $M$ is
compact and $i$ is an injective immersion, then $i$ is an
embedding. This is not the case if $M$ is not compact since $i(M)$
can be self-approaching.
\begin{definition}
A \textbf{Riemannian manifold} $(M,g)$ is a smooth manifold $M$
together with a tensor\footnote{ A tensor $T$ of type $(m,n)$ is a
multilinear form $T_p M\times \ldots T_p M \times T^*_p M \times
\ldots \times T^*_p M \rightarrow \R$ ($n$-times $T_p M$,
$m$-times $T^*_p M$).} of type $(0,2)$, called the metric tensor
$g$, at each $p \in M$, such that $g$ defines an inner product on
the tangent space $T_p M$ which varies smoothly over $M$. The
volume form induced by $g$ is given in local coordinates as
$dV=\sqrt{\det g} \, dx^1 \wedge \ldots \wedge dx^m$. $dV$ is
uniquely determined by $dV(e_1,\ldots,e_m)=1$ for any oriented
orthonormal basis $e_1,\ldots,e_m$ in $T_x M$.
\end{definition}
The metric tensor induces for every $p \in M$ an isometric
isomorphism between the tangent space $T_p M$ and its dual $T^*_p
M$. A submanifold $M$ of a Riemannian manifold $(X,g)$ has a
natural Riemannian metric $h$ induced from $X$ in the following
way. Let $i:M \rightarrow X$ be an embedding so that $M$ is a
submanifold of $X$. Then one can induce a metric $h$ on $M$ using
the mapping $i$, namely $h=i^*g$, where $i^*:T^*_{i(x)}X
\rightarrow T^*_xM$ is the pull-back\footnote{$T^*_x M$
 is the dual of the tangent space $T_x M$. Every differentiable mapping $i:M\rightarrow X$
induces a pull-back $i^*:T^*_{i(x)}X\rightarrow T^*_x M$. Let $u
\in T_x M$, $w \in T^*_{i(x)} X$ and denote by $i'$ the
differential of $i$. Then $i^*$ is defined by $(i^*w)(u)=w(i'u)$.}
of the differentiable mapping $i$. In this case $i$ trivially is
an isometric embedding of $(M,h)$ into $(X,g)$. In the paper we
always use on the submanifold $M$ the metric induced from $\R^d$.
\begin{definition}
The \textbf{Laplace-Beltrami operator $\Delta_M$} of a Riemannian
manifold is defined as $\Delta_M=\div(\grad)$. For a twice
differentiable function $f: M \rightarrow \R$ it is explicitly
given as
\[ \Delta_M f=\frac{1}{\sqrt{\det g}}\frac{\partial}{\partial x^j}\Big(\sqrt{\det g}\,g^{ij}\frac{\partial f}{\partial x^i}\Big), \]
where $g^{ij}$ are the components of the inverse of the metric
tensor $g=g_{ij}\,dx^i \otimes dx^j$.
\end{definition}

\subsection{Normal Coordinates}
Since in the proofs we use normal coordinates, we give here a
short introduction. Intuitively, normal coordinates around a point
$p$ of an $m$-dimensional Riemannian manifold $M$ are coordinates
chosen such that $M$ looks around $p$ like $\R^m$ in the best
possible way. This is achieved by adapting the coordinate lines to
geodesics through the point $p$. The reference for the following
material is the book of \cite{Jost2002}. We denote by $c_v$ the
unique geodesic starting at $c(0)=x$ with tangent vector
$\dot{c}(0)=v$ ($c_v$ depends smoothly on $p$ and $v$).
\begin{definition}
Let $M$ be a Riemannian manifold, $p \in M$, and $V_p= \{ v \in
T_p M,\; c_v \;\textrm{defined}$ $\; \textrm{on} \; [0,1]\}$,
then,
 $\exp_p: V_p \rightarrow M$, $v \mapsto c_v(1)$,
is called the \textbf{exponential map} of $M$ at $p$.
\end{definition}
It can be shown that $\exp_p$ maps a neighborhood of $0 \in T_p M$
diffeomorphically onto a neighborhood $U$ of $p \in M$. This
justifies the definition of normal coordinates.
\begin{definition}
Let $U$ be a neighborhood of $p$ in $M$ such that $\exp_p$ is a
diffeomorphism. The local coordinates defined by the chart
$(U,\exp_p^{-1})$ are called \textbf{normal coordinates} at $p$.
\end{definition}
Note that in $T_p M\simeq \R^m \supset \exp_p^{-1}(U)$ we use
always an orthonormal basis. The injectivity radius describes the
largest ball around $p$ such that normal coordinates can be
introduced.
\begin{definition}
Let $M$ be a Riemannian manifold. The \textbf{injectivity radius}
of $p \in M$ is
\[ \inj(p)= \sup\{ \rho > 0,\; \exp_p \; \mathrm{is \; defined\; on} \; \overline{B_{\R^m}(0,\rho)} \; \mathrm{and \; injective}\}.\]
\end{definition}
It can be shown that $\inj(p)>0, \forall p \in M \backslash
\partial M$. Moreover,  for compact manifolds without boundary
there exists a lower bound $\inj_{\min}>0$ such that $\inj(p)\geq
\inj_{\min}, \forall p \in M$. However,  for manifolds with
boundary one has $\inj(p_n) \rightarrow 0$ for any sequence of
points $p_n$ with limit on the boundary. The motivation for
introducing normal coordinates is that the geometry is
particularly simple in these coordinates. The following theorem
makes this more precise.
\begin{theorem}
In normal coordinates around $p$ one has for the Riemannian metric
$g$
and the Laplace-Beltrami operator $\Delta_M$ applied to a function
$f$ at $p=\exp^{-1}_p(0)$,
\[ g_{ij}(0)=\delta_{ij}, \quad \frac{\partial}{\partial x^k}g_{ij}(0)=0, \quad (\Delta_M f)(0)= \sum_{i=1}^m \frac{\partial^2 f}{\partial (x^i)^2}(0)\] 
\end{theorem}
The second derivatives of the metric tensor cannot be made to
vanish in general. There curvature effects come into play which
cannot be deleted by a coordinate transformation. To summarize,
normal coordinates with center $p$ achieve that, up to first
order, the geometry of $M$ at point $p$ looks like that of $\R^m$.

\subsection{The Second Fundamental Form}\label{sec:second-fundamental-form}
In this section we assume that $M$ is an isometrically embedded
submanifold of a manifold $X$. At each point $p \in M$ one can
decompose the tangent space $T_pX$ into a subspace $T_p M$, which
is the tangent space to $M$, and the orthogonal normal space $N_p
M$. In the same way one can split the covariant derivative of $X$
at $p$, $\tilde{\nabla}_U V$ into a component tangent
$(\tilde{\nabla}_U V)^\top$ and normal $(\tilde{\nabla}_U V)^\bot$
to $M$.
\begin{definition}
The \textbf{second fundamental form} $\Pi$ of an isometrically
embedded submanifold $M$ of $X$ is defined as
\[ \Pi: T_p M \otimes T_p M \rightarrow N_p M, \quad \Pi(U,V)=(\tilde{\nabla}_U V)^\bot \]
\end{definition}
The following theorem, see \cite{Lee1997}, then shows that the
covariant derivative of $M$ at $p$ is nothing else than the
projection of the covariant derivative of $X$ at $p$ onto $T_p M$.
\begin{theorem}[Gauss Formula] \label{th:gauss-formula}
Let $U,V$ be vector fields on $M$ which are arbitrarily extended
to $X$, then the following holds along $M$
\[ \tilde{\nabla}_U V =\nabla_U V + \Pi(U,V) \]
where $\tilde{\nabla}$ is the covariant derivative of $X$ and
$\nabla$ the covariant derivative of $M$.
\end{theorem}
The second fundamental form connects also the curvature tensors of
$X$ and $M$.
\begin{theorem}[Gauss equation] \label{th:gauss-equation}
For any $U,V,W,Z \in T_p M$ the following equation holds
\[ \tilde{R}(U,V,W,Z)=R(U,V,W,Z)-\inner{\Pi(U,Z),\Pi(V,W)} + \inner{\Pi(U,W),\Pi(V,Z)}, \]
where $\tilde{R}$ and $R$ are the Riemann curvature\footnote{The
Riemann curvature tensor of a Riemannian manifold $M$ is defined
as $R:T_p M \otimes T_p M \otimes T_p M \rightarrow T^*_p M$,
\[ R(X,Y)Z = \nabla_X \nabla_Y  Z - \nabla_Y \nabla_X Z -\nabla_{[X,Y]} Z. \]
In local coordinates $x^i$, $R_{ijk}\,^l \partial_l =
R(\partial_i, \partial_j)\partial_k$ and
$R_{ijkm}=g_{lm}R_{ijk}\,^l$.} tensors of $X$ and $M$.
\end{theorem}
In this paper we derive a relationship between distances in $M$
and the corresponding distances in $X$. Since Riemannian manifolds
are length spaces and therefore the distance is induced by length
minimizing curves (locally the geodesics), it is of special
interest to connect properties of curves of $M$ with respect to
$X$. Applying the Gauss Formula to a curve $c(t):(t_0,t_1)
\rightarrow M$ yields the following
\[ \tilde{D}_t V = D_t V + \Pi(V,\dot{c}), \]
where $\tilde{D}_t = \dot{c}^a \tilde{\nabla}_a$ and $\dot{c}$ is
the tangent vector field to the curve $c(t)$. Now let $c(t)$ be a
geodesic parameterized by arc-length, that is with unit-speed,
then its acceleration fulfills $D_t \dot{c}=\dot{c}^a \nabla_a
\dot{c}^b=0$ (however that is only true locally in the interior of
$M$, globally if $M$ has boundary length minimizing curves may
behave differently especially if a length minimizing curve goes
along the boundary its acceleration can be non-zero), and one gets
for the acceleration in the ambient space
\[ \tilde{D}_t \dot{c} =  \Pi(\dot{c},\dot{c}). \]
In our setting where $X=\R^d$ the term $\tilde{D}_t \dot{c}$ is
just the ordinary acceleration $\ddot{c}$ in $\R^d$. Remember that
the norm of the acceleration vector is inverse to the curvature of
the curve at that point (if $c$ is parameterized by
arc-length\footnote{Note that if $c$ is parameterized by
arc-length, $\dot{c}$ is tangent to $M$, that is in particular
$\norm{\dot{c}}_{T_x X}=\norm{\dot{c}}_{T_x M}$}). Due to this
connection it becomes more apparent why the second fundamental
form is often called the extrinsic curvature (with respect to $X$).\\
The following Lemma shows that the second fundamental form $\Pi$
of an isometrically embedded submanifold $M$ of $\R^d$ is in
normal coordinates just the Hessian of $i$.
\begin{lemma}\label{le:sform-normal}
Let $e_{\alpha},\, \alpha=1,\ldots,d$ denote an orthonormal basis
of $T_{i(x)}\R^d$ then the second fundamental form of $M$ in
normal coordinates $y$ is given as:
\[ \Pi(\partial_{y^i},\partial_{y^j})\Big|_{0}= \frac{\partial^2 i^{\alpha}}{\partial y^i \partial y^j}e_{\alpha}.\]
\end{lemma}
\begin{proof}
Let $\tilde{\nabla}$ be the flat connection of $\R^d$ and $\nabla$
the connection of $M$. Then by Theorem \ref{th:gauss-formula},
$\Pi(\partial_{y^i},\partial_{y^j})=\tilde{\nabla}_{i^*\partial_{y^i}}(i^*\partial_{y^j})-\nabla_{\partial_{y^i}}\partial_{y^j}
=
\partial_{y^i}\left(\frac{\partial i^\alpha }{\partial_{y^j}}\right)e_{\alpha}=\frac{\partial^2 i^{\alpha}}{\partial y^i \partial y^j}e_{\alpha},$
where the second equality follows from the flatness of
$\tilde{\nabla}$ and $\Gamma^i_{\,jk}\Big|_{0}=0$ in normal
coordinates.
\end{proof}

\section{Proofs and Lemmas}\label{sec:proofs}
\subsection*{Proof of Proposition \ref{pro:averaging-op-manifold}}
The following lemmas are needed in the proof.
\begin{lemma}\label{le:partial-integration}
If the kernel $k:\R_+
\rightarrow \R_+$ satisfies Assumptions \ref{kernel-assumptions},
then
\begin{equation}
\int_{\R^m}\frac{\partial k}{\partial x}(\norm{u}^2)u^i u^j u^k
u^l du
   =-\frac{1}{2}C_2 \left[\delta^{ij}\delta^{kl}+\delta^{ik}\delta^{jl}+\delta^{il}\delta^{jk}\right].
\end{equation}
\end{lemma}
\begin{proof}
Note first that for a function $f(\norm{u}^2)$ one has
$\frac{\partial f}{\partial \norm{u}^2}=\frac{\partial f}{\partial
u^2_i}$. The rest follows from partial integration.
\begin{align*}
  \int\limits_{-\infty}^\infty \frac{\partial k}{\partial u^2}(u^2) \,u^2\, du
    = \int\limits_{0}^\infty \frac{\partial k}{\partial v}(v)\, \sqrt{v}\, dv
    = \left[k(v)\,\sqrt{v}\right]_0^\infty - \int\limits_{0}^\infty k(v)\, \frac{1}{2\sqrt{v}} \,dv
    =-\frac{1}{2}\int\limits_{-\infty}^\infty k(u^2) du,
\end{align*}
where $\left[k(v)\,\sqrt{v}\right]_0^\infty=0$ due to the boundedness and exponential decay of $k$.\\
In the same way one can derive, $\int_{-\infty}^\infty
\frac{\partial k}{\partial u^2}(u^2) \,u^4\,
du=-\frac{3}{2}\int_{-\infty}^\infty k(u^2)u^2\, du$. The result
follows by noting that since $k$ is an even function only
integration over even powers of coordinates will be non-zero.
\end{proof}
\begin{lemma}\label{le:averaging}
Let $k$ satisfy Assumption \ref{kernel-assumptions} and let
$V_{ijkl}$ be a given tensor. Assume now $\norm{z}^2\geq
\norm{z}^2 + V_{ijkl}z^i z^j z^k z^l + \beta(z) \norm{z}^5 \geq
\frac{1}{4}\norm{z}^2$ on $B(0,r_\mathrm{min}) \subset \R^m$,
where $\beta(z)$ is con\-ti\-nuous and $\beta(z) \sim O(1)$ as
$z\rightarrow 0$. Then there exists a constant $C$ and a $h_0>0$
such that for all $h<h_0$ and all $f \in
C^3(B(0,r_\mathrm{min}))$,
\begin{align*}
&\Big|\int_{B(0,r_\mathrm{min})} k_h\left(\frac{\norm{z}^2+V_{ijkl}z^i z^j z^k z^l + \beta(z) \norm{z}^5)}{h^2}\right)f(z) dz \nonumber \\
&- \Big( C_1 f(0) + C_2 \frac{h^2}{2}\Big[(\Delta f)(0) -
f(0)\sum_{i,k}^m V_{iikk}+V_{ikik}+V_{ikki}\Big]\Big)\Big|\leq C
h^3.
\end{align*}
where $C$ is a constant depending on $k$, $r_\mathrm{min}$,
$V_{ijkl}$ and $\norm{f}_{C^3}$.
\end{lemma}
\begin{proof}
As a first step we do a Taylor expansion of the kernel around
$\norm{z}^2/h^2$:
\begin{align*}
 k_h\bigg(\frac{\norm{z}^2 + \eta}{h^2}\bigg)
=k_h\bigg(\frac{\norm{z}^2}{h^2}\bigg) + \frac{\partial
k_h}{\partial
x}\Big|_{\frac{\norm{z}^2}{h^2}}\frac{\eta}{h^2}+\frac{\partial^2
  k_h(x)}{\partial x^2}\Big|_{\frac{\norm{z}^2(1-\theta)+\theta\, \eta}{h^2}}\frac{\eta^2}{h^4},
\end{align*}
where in the last term $0\leq \theta(z) \leq 1$. We then decompose
the integral:
\begin{align*}
&\int_{B(0,r_\mathrm{min})} k_h\bigg(\frac{\norm{z}^2+V_{ijkl}z^i z^j z^k z^l + \beta(z)\norm{z}^5}{h^2}\bigg)f(z) dz\\
=& \hspace{-1.5mm}\int\limits_{\R^m}
\hspace{-1.5mm}\bigg(k_h\bigg(\frac{\norm{z}^2}{h^2}\bigg) +
  \frac{\partial k_h}{\partial x}\Big|_{\frac{\norm{z}^2}{h^2}}\frac{V_{ijkl}\,z^i z^j z^k z^l}{h^2}\bigg)
\Big(f(0)+\inner{\nabla f\big|_0,z} + \frac{1}{2} \frac{\partial^2
f}{\partial z^i \partial z^j}\Big|_0 z^i z^j\Big)d z
   + \sum_{i=0}^4 \alpha_i,
\end{align*}
where we define the five error terms $\alpha_i$ as:
\begin{align*}
\alpha_0 &= \int_{B(0,r_\mathrm{min})} \frac{\partial k_h}{\partial x}\Big|_{\frac{\norm{z}^2}{h^2}} \frac{\beta(z)\norm{z}^5}{h^2}f(z)\,dz, \\
\alpha_1 &= \int_{B(0,r_\mathrm{min})} \frac{\partial^2
k_h}{\partial x^2}\Big|_{\frac{\norm{z}^2(1-\theta) + \theta
\eta}{h^2}}
            \frac{\Big(V_{ijkl}z^i z^j z^k z^l +\beta(z)\norm{z}^5 \Big)^2}{h^4}\, f(z)dz,\\
\alpha_2&=
\int_{B(0,r_\mathrm{min})}k_h\bigg(\frac{\norm{z}^2+V_{ijkl}z^i
z^j z^k z^l + \beta(z) \norm{z}^5}{h^2}\bigg)
             \frac{1}{6}\frac{\partial^3 f}{\partial z^i \partial z^j \partial z^k}(\theta z) z^i z^j z^k dz,\\
\alpha_3 &=  \int_{\R^m\backslash B(0,r_\mathrm{min})}
k_h\bigg(\frac{\norm{z}^2}{h^2}\bigg)
             \bigg(f(0)+\inner{\nabla f\big|_0,z} + \frac{1}{2} \frac{\partial^2 f}{\partial z^i \partial z^j}\Big|_0z^i z^j\bigg)dz,
\end{align*}
\begin{align*}
\alpha_4 &=  \int_{\R^m\backslash B(0,r_\mathrm{min})}
              \frac{\partial k_h}{\partial x}\Big|_{\frac{\norm{z}^2}{h^2}}\frac{V_{ijkl}\,z^i z^j z^k z^l}{h^2}
              \bigg(f(0)+\inner{\nabla f\big|_0,z}
              + \frac{1}{2}\frac{\partial^2 f}{\partial z^i \partial z^j}\Big|_0 z^i z^j\bigg)dz,
\end{align*}
where in $\alpha_1$, $\eta=V_{ijkl}z^i z^j z^k z^l + \beta(z) \norm{z}^5$. With $\int_{\R^m}k(\norm{z}^2)\,z_i \,dz =0, \, \forall \,i$, and\\
$\int_{\R^m} k(\norm{z}^2)\,z_i\, z_j dz=0$ if $i\neq j$, and
Lemma \ref{le:partial-integration} the main term simplifies to:
\begin{align*}
   & \int_{\R^m} \bigg(k_h\bigg(\frac{\norm{z}^2}{h^2}\bigg) +\frac{\partial k_h(x)}{\partial x}\Big|_{\frac{\norm{z}^2}{h^2}}
     \frac{V_{ijkl}\,z^i z^j z^k z^l}{h^2}\bigg)
     \bigg(f(0)+ \frac{1}{2} \frac{\partial^2 f}{\partial z^i \partial z^j }\Big|_0 z^i z^j \bigg) dz\\
  =& \int_{\R^m} \Big(k\big(\norm{u}^2\big) + h^2 \frac{\partial k(x)}{\partial x}\Big|_{\norm{u}^2}V_{ijkl}\,u^i u^j u^k u^l\Big)
     \Big(f(0)+ \frac{h^2}{2} \frac{\partial^2 f}{\partial z^i \partial z^j }\Big|_0 u^i u^j \Big) du\\
  =&\, C_1 f(0) - \frac{h^2}{2}C_2 f(0) V_{ijkl}\Big[\delta^{ij}\delta^{kl}+\delta^{ik}\delta^{jl}+\delta^{il}\delta^{jk}\Big]
     + \frac{h^2}{2} C_2 \sum_{i=1}^m \frac{\partial^2 f}{\partial (z^i)^2}\Big|_0 + O(h^4)
\end{align*}
where the $O(h^4)$ term is finite due to the exponential decay of
$k$ and depends on $k$, $r_{\min}$, $V_{ijkl}$ and
$\norm{f}_{C^3}$. Now we can upper bound the remaining error terms
$\alpha_i,\, i=0,\ldots,4$. For the argument of the kernel in
$\alpha_1$ and $\alpha_2$ we have by our assumptions
 on $B(0,r_\mathrm{min})$:
\[ \frac{\norm{z}^2}{h^2}\geq\frac{\norm{z}^2+V_{ijkl}z^i z^j z^k z^l + \beta(z)\norm{z}^5}{h^2} \geq \frac{\norm{z}^2}{4h^2}.\]
Note that this inequality implies that $\beta$ is uniformly
bounded on $B(0,r_{\min})$ in terms of $r_{\min}$ and $V_{ijkl}$.
Moreover,  for small enough $h$ we have
$\frac{r_\mathrm{min}}{h}\geq \sqrt{A}$ (see Assumptions
\ref{kernel-assumptions} for the definition of $A$) so that we can
use the exponential decay of $k$ for $\alpha_3$ and $\alpha_4$.
\begin{align*}
|\alpha_0| &\leq h^3
\norm{f}_{C^3}\int_{B(0,\frac{r_\mathrm{min}}{h})} \frac{\partial
k_h}{\partial x}\bigg|_{\norm{u}^2} |\beta(hu)|\norm{u}^5\,du
\end{align*}
Since $\frac{\partial k_h}{\partial x}$ is bounded and has
exponential decay, one has $|\alpha_0|\leq K_0 \,h^3$ where $K_0$
depends on $k$, $r_{\min}$ and $\norm{f}_{C^3}$.
\begin{align*}
|\alpha_1|&\leq \int_{B(0,r_\mathrm{min})}\bigg| \frac{\partial^2
k_h}{\partial x^2}\bigg(\frac{\norm{z}^2(1-\theta) + \theta
\eta}{h^2}\bigg)\bigg| \frac{\big(V_{ijkl}z^i z^j z^k z^l +\beta(z)\norm{z}^5 \big)^2}{h^4} f(z)dz \\
          \leq & h^4 \norm{f}_{C^3}    \int_{B(0,\frac{r_{\min}}{h})}  \bigg|
          \frac{\partial^2 k}{\partial x^2}\big(\norm{u}^2(1-\theta) + \theta \eta\big)\bigg|
          \Big(m^2 \max_{i,j,k,l}|V_{ijkl}|
          \norm{u}^4 +h \norm{\beta}_\infty \norm{u}^5 \Big)^2 du
\end{align*}
First suppose $\frac{r_{\mathrm{\min}}}{h}\leq 2\sqrt{A}$ then the
integral is bounded since the integrands are bounded on
$B(0,\frac{r_{\min}}{h})$. Now suppose
$\frac{r_{\mathrm{\min}}}{h}\geq 2\sqrt{A}$ and decompose
$B(0,\frac{r_\mathrm{min}}{h})$ as
$B(0,\frac{r_{\mathrm{min}}}{h}) = B(0,2\sqrt{A}) \cup
B(0,\frac{r_\mathrm{min}}{h}) \backslash B(0,2\sqrt{A})$. On
$B(0,2\sqrt{A})$ the integral is finite since
$\big|\frac{\partial^2 k}{\partial x^2}\big|$ is bounded and on
the complement the integral is also finite since
$\big|\frac{\partial^2 k}{\partial x^2}\big|$ has exponential
decay since by assumption
\[ \norm{u}^2(1-\theta(hu)) + \theta(hu) \eta(hu) \geq \frac{1}{4}\norm{u}^2 \geq A.\]
Therefore there exists a constant $K_1$ such that $|\alpha_1|\leq
K_1 \,h^4$.
\begin{align*}
|\alpha_2|&\leq
\hspace{-6mm}\int\limits_{B(0,r_\mathrm{min})}\hspace{-5mm}k_h\bigg(\frac{\norm{z}^2}{4h^2}\bigg)
\frac{1}{6}\frac{\partial^3 f}{\partial z^i \partial z^j
             \partial z^k}(\theta z) z^i z^j z^k dz
          \leq \frac{m^{3/2}\norm{f}_{C^3} h^3}{6} \int\limits_{\R^m} \hspace{-2mm} k\bigg(\frac{\norm{u}^2}{4}\bigg)\norm{u}^3 du \leq K_2\, h^3,
\end{align*}
\begin{align*}
|\alpha_3|&\leq \int_{\R^m\backslash B(0,r_\mathrm{min})}
k_h\bigg(\frac{\norm{z}^2}{h^2}\bigg)
             \Big(f(0)+\inner{\nabla f\big|_0,z} + \frac{1}{2} \frac{\partial^2 f}{\partial z^i \partial z^j}\Big|_0 z^i z^j\Big)\,dz\\
          &\leq c \norm{f}_{C^3}\hspace{-2mm} \int\limits_{\R^m\backslash B(0,r_\mathrm{min})}\hspace{-2mm}e^{-\frac{\alpha \norm{z}^2}{h^2}}(1+m\, h^2\norm{z}^2)\, dz
          \leq c\, e^{-\alpha\frac{r^2_\mathrm{min}}{2h^2}}\bigg(\frac{2\pi}{\alpha}\bigg)^\frac{m}{2}
            \Big(1+m\, h^2\frac{m}{\alpha}\Big),\\
|\alpha_4|&\leq K_4 \hspace{-7mm}\int\limits_{\R^m\backslash
B(0,r_\mathrm{min})} \hspace{-5mm} \Big|\frac{\partial
k_h}{\partial x}\Big(\frac{\norm{z}^2}{h^2}\Big)\Big|
            \frac{\norm{z}^4+\norm{z}^6}{h^2}\,dz
          \leq c \, K_4\, h^2 e^{-\alpha\frac{r^2_\mathrm{min}}{2h^2}}
            \int\limits_{\R^m} \hspace{-1mm} e^{-\alpha\norm{u}^2} (\norm{u}^4 + h^2\norm{u}^6)du,
\end{align*}
where $K_4$ is a constant depending on $\max_{i,j,k,l}|V_{ijkl}|$
and $\norm{f}_{C^3}$. Now one has\footnote{This inequality can be
deduced from $e^x\geq x^n$ for all $x\geq 4n^2$. }:
$e^{-\frac{\xi^2}{h^2}}\leq h^s/\xi^s$ for $h\leq \xi/s$. in
particular,  it holds $h^3\geq
e^{-\alpha\frac{r^2_\mathrm{min}}{2h^2}}$ for $h\leq
\frac{1}{3}\sqrt{\frac{\alpha}{2}}\, r_{\min}$, so that for
$h<\min\{\frac{1}{3}\sqrt{\frac{\alpha}{2}}\, r_{\min},
\frac{r_{\min}}{\sqrt{A}}\}=h_0$ all error terms are smaller than
a constant times $h^3$ where the constant depends on $k$,
$r_\mathrm{min}$, $V_{ijkl}$ and $\norm{f}_{C^3}$. This finishes
the proof.
\end{proof}
Now we are ready to prove Proposition
\ref{pro:averaging-op-manifold},
\begin{proof}
Let
$\epsilon=\frac{1}{3}\min\{\mathrm{inj}(x),\pi\rho\}$\footnote{The
factor $1/3$ is needed in Theorem
\ref{th:convergence-bias-manifold}} where $\epsilon$ is positive
by the assumptions on $M$. Then we decompose $M$ as
$M=B(x,\epsilon) \cup (M \backslash B(x,\epsilon))$ and integrate
separately. The integral over $M \backslash B(x,\epsilon)$ can be
upper bounded by using the definition of $\delta(x)$ (see
Assumption \ref{manifold-assumptions}) and the fact that $k$ is
non-increasing:
\begin{align*}\label{eq:manifold-decomposition}
  \int\limits_{M}k_h\big(\norm{i(x)-i(y)}^2_{\R^d}\big)f(y)p(y)\sqrt{\det g}\,dy
=\hspace{-3mm}\int\limits_{B(x,\epsilon)}k_h\big(\norm{i(x)-i(y)}^2_{\R^d}\big)f(y)p(y)\sqrt{\det g}\,dy  \nonumber \\
 + \int\limits_{M \backslash B(x,\epsilon)}k_h\big(\norm{i(x)-i(y)}^2_{\R^d}\big)f(y)p(y)\sqrt{\det g}\,dy
\end{align*}
Since $k$ is non-increasing, we have the following inequality for
the integral over $M\backslash B(x,\epsilon)$:
\begin{align*}
&  \int_{M \backslash
B(x,\epsilon)}k_h\Big(\norm{i(x)-i(y)}^2_{\R^d}\Big)f(y)p(y)\sqrt{\det
g}\,dy\leq
\frac{1}{h^m}k\bigg(\frac{\delta(x)^2}{h^2}\bigg)\norm{f}_{\infty}
\nonumber
\end{align*}
Since $\delta(x)$ is positive by assumption and $k$ decays
exponentially, we can make the upper bound smaller than $h^3$ for
small enough $h$. Now we deal with the integral over
$B(x,\epsilon)$. Since $\epsilon$ is smaller than the injectivity
radius $\mathrm{inj}(x)$, we can introduce normal coordinates
$z=\exp_x^{-1}(y)$ on $B(x,\epsilon)$, so that we can rewrite the
integral using Proposition \ref{pro:normal-coordinates} as:
\begin{equation}\label{eq:proof1-1}
\int_{B(0,\epsilon)}k_h\bigg(\frac{\norm{z}^2-\frac{1}{12}\sum_{\alpha=1}^d
\frac{\partial^2 i^\alpha}{\partial z^a
\partial z^b}\frac{\partial^2 i^\alpha}{\partial z^u \partial z^v}z^a z^b z^u z^v + O(\norm{z}^5)}{h^2}\bigg)p(z)f(z)\sqrt{\det g}\,dz
\end{equation}
Using our assumptions, we see that  $pf\sqrt{\det g}$ is in
$C^3(B(0,\epsilon))$. Moreover, by Corollary
\ref{co:comparison-ex-in} one has for $d_M(x,y)\leq \pi \rho$,
$\frac{1}{2}d_M(x,y)\leq \norm{x-y} \leq d_M(x,y)$. Therefore we
can apply Lemma \ref{le:averaging} and compute the integral in
\eqref{eq:proof1-1} which results in:
\begin{align}\label{eq:proof1-2}
 p(0)f(0) \Big(C_1 +\frac{h^2\,C_2}{24}\sum_{\alpha=1}^d
  \frac{\partial^2 i^\alpha}{\partial z^a \partial z^b}\frac{\partial^2 i^\alpha}{\partial z^c \partial z^d}
  \Big[\delta^{ab}\delta^{cd}+\delta^{ac}\delta^{bd}+\delta^{ad}\delta^{bc}\Big]\Big) \nonumber \\
+ \frac{h^2\,C_2}{2} \Delta_M(pf\sqrt{\det g})\Big|_{0} + O(h^3),
\end{align}
where we have used that in normal coordinates $z^i$ at $0$ the
Laplace-Beltrami operator $\Delta_M$ is given as $\Delta_M
f\Big|_x=\sum_{i=1}^m \frac{\partial^2 f}{\partial
(z^i)^2}\Big|_0$. The second term in the above equation can be
evaluated using the Gauss equations, see \cite[Proposition
6]{SmoWeiWit2004}.
\begin{align*}
\sum_{a,b=1}^m& \sum_{\alpha=1}^d \frac{\partial^2
i^\alpha}{\partial z^a \partial z^b}\frac{\partial^2
i^\alpha}{\partial z^c \partial
z^d}\Big[\delta^{ab}\delta^{cd}+\delta^{ac}\delta^{bd}+\delta^{ad}\delta^{bc}\Big]
\hspace{-1.8mm}=\hspace{-1.8mm}\sum_{a,b=1}^m \sum_{\alpha=1}^d
\frac{\partial^2 i^\alpha}{\partial (z^a)^2}\frac{\partial^2
i^\alpha}{\partial (z^b)^2} + 2 \frac{\partial^2
i^\alpha}{\partial z^a
\partial z^b}\frac{\partial^2 i^\alpha}{\partial z^a \partial z^b} \\
=& 2 \sum_{a,b=1}^m \sum_{\alpha=1}^d \left(\frac{\partial^2
i^\alpha}{\partial z^a
\partial z^b}\frac{\partial^2 i^\alpha}{\partial z^a \partial
z^b} - \frac{\partial^2 i^\alpha}{\partial
(z^a)^2}\frac{\partial^2 i^\alpha}{\partial (z^b)^2}\right) + 3
\sum_{a,b=1}^m \sum_{\alpha=1}^d
\frac{\partial^2 i^\alpha}{\partial (z^a)^2}\frac{\partial^2 i^\alpha}{\partial (z^b)^2}\\
=& 2 \sum_{a,b=1}^m
\inner{\Pi(\partial_{z^a},\partial_{z^b}),\Pi(\partial_{z^a},\partial_{z^b})}-
\inner{\Pi(\partial_{z^a},\partial_{z^a}),\Pi(\partial_{z^b},\partial_{z^b})} + 3 \bigg \|\sum_{a=1}^m \Pi(\partial_{z^a},\partial_{z^a})\bigg \|_{T_{i(x)}\R^d}^2\\
=& - 2 R + 3 \bigg \|\sum_{j=1}^m
\Pi(\partial_{z^j},\partial_{z^j}) \bigg\|_{T_{i(x)}\R^d}^2,
\nonumber
\end{align*}
where $R$ is the scalar curvature and we used Lemma
\ref{le:sform-normal} in the third equality. Plugging this result
into \eqref{eq:proof1-2} and using from Proposition
\ref{pro:normal-coordinates}, $\Delta_M \sqrt{\det g}\big|_0
=-\frac{1}{3}R$, we are done.
\end{proof}
\label{app:theorem}
\begin{lemma}\label{le:averaged-density-bounds}
Let $k$ have compact support on $[0,R_k^2]$ and let $0<h \leq
h_{\max}$. Then for any $x \in M$ there exist constants $D_1,
D_2>0$ independent of $h$ such that for any $y \in
B_{\R^d}(x,hR_k) \cap M$,
\[  0<D_1 \leq p_h(y) \leq D_2.  \]
\end{lemma}
\begin{proof}
First suppose that $hR_k < s\bydef \min\{\kappa/2,R_0/2\}$. Since
$\norm{y-z}\leq hR_k\leq \kappa/2$ we have by Lemma
\ref{le:comparison-ex-in}: $\frac{1}{2}d_M(y,z)\leq \norm{y-z}\leq
d_M(y,z)$. Moreover,  since $p(x)>0$ on $M$ and $p$ is bounded and
continuous, there exist lower and upper bounds $p_{\min}$ and
$p_{\max}$ on the density on $B_M(x,4hR_k)$. That implies
\begin{align*}
p_h(y)\leq \frac{\norm{k}_\infty}{h^m} p_{\max}
\int_{B_M(y,2hR_k)} \sqrt{\det g}\, dz \leq \norm{k}_\infty
p_{\max}\, S_2 \, 2^m R_k^m,
\end{align*}
where the last inequality follows from Lemma
\ref{le:volume-bounds}. Note further that $d_M(x,y)\leq 2hR_k$ and
$d_M(y,z)\leq 2hR_k$ implies $d_M(x,z)\leq 4hR_k$. Since the
kernel function is continuous there exists an $r_k$ such that
$k(x)\geq \norm{k}_\infty/2$ for $0<x\leq r_k$. We get
\begin{align*}
p_h(y)&\geq \frac{\norm{k}_\infty}{2h^m} \int\limits_{B_{\R^d}(x,h
\,r_k) \cap M} \hspace{-5mm} p(z)\sqrt{\det g}\, dz \geq
\frac{\norm{k}_\infty}{2h^m} p_{\min} \vol_M(B_M(x,h\,r_k))\geq
\frac{\norm{k}_\infty}{2}p_{\min}\,S_1 \,r_k^m.
\end{align*}
Now suppose $s\leq hR_k$ and $h\leq h_{\max}$. Then $p_h(y) \leq
\frac{\norm{k}_\infty}{h^m} \leq  \norm{k}_\infty
\left(\frac{R_k}{s}\right)^m$. For the lower bound we get
\begin{align*} p_h(y) &\geq \int_M k_h(d_M(y,z)) p(z)\sqrt{\det g}\, dz \geq \int_{B_M(y,h\,r_k)}k_h(d_M(y,z))p(z)\sqrt{\det g}\,dz\\
&\geq \frac{\norm{k}_\infty}{2h^m}\,\Pr\Big(B_M(y,h\,r_k)\Big)\geq
\frac{\norm{k}_\infty}{2h_{\max}^m}\,
\Pr\Big(B_M(y,s\,\frac{r_k}{R_k})\Big)
\end{align*}
Since $p$ is continuous and $p>0$, the function $y \rightarrow
\Pr\Big(B_M(y,s\,\frac{r_k}{R_k})\Big)$ is continuous and positive
and therefore has a lower bound greater zero on the ball
$B_{\R^d}(x,hR_k) \cap M$.
\end{proof}
\end{appendix}

\vskip 0.2in
\bibliographystyle{apalike}
\bibliography{../../regul}
\end{document}